\documentclass[10pt]{article}
\usepackage{amssymb,graphicx, cancel,ulem}
\usepackage{amsmath, amscd,amsthm, tikz-cd,color,amsfonts, amssymb,verbatim,caption}

\definecolor{grey}{rgb}{0.6,0.6,0.6}
\definecolor{darkgreen}{rgb}{0.2,0.5,0.0}
\definecolor{darkblue}{rgb}{0.2,0.2,0.71}

  \newenvironment{packed_enum}{
\begin{enumerate}
  \setlength{\itemsep}{0pt}
  \setlength{\parskip}{0pt}
  \setlength{\parsep}{0pt}
}{\end{enumerate}}

\usepackage[colorlinks=true, pdfstartview=FitV, linkcolor=darkblue, citecolor=darkblue,urlcolor=darkblue]{hyperref}

\def\Tm{T_-}
\def\Tp{T_+}
\def\T3{T_3}

\def \CC{\mathcal C}
\def \SS{\mathfrak S}
\newlength{\dinwidth}
\def\le{\left}
\def \QED{\hfill $\Box$\par \vskip 4pt}
\def\ri{\right}
\def \P{\mathcal P}
\newlength{\dinmargin}
\def\bea{\begin{eqnarray}}
\def\eea{\end{eqnarray}}   
\def \wt{ \widetilde }
\def \&{\hspace{-15pt}&}

\def \d{{\mathrm d}}
\def\res{\mathop{ \mathrm {res}}}

\newtheorem{definition}{Definition}[section]
\newtheorem{theorem}{Theorem}[section]
\newtheorem{proposition}{Proposition}[section]
\newtheorem{corollary}{Corollary}[section]
\newtheorem{remark}{Remark}[section]
\newtheorem{lemma}{Lemma}[section]

\newtheorem{conjecture}{Conjecture}

\def\be{\begin{equation}}
\def\ee{\end{equation}}
\def\ben{\begin{displaymath}}
\def\een{\end{displaymath}}
\def\baa{\begin{eqnarray}}
\def\eaa{\end{eqnarray}}

\def\ba{\begin{array}}
\def\ea{\end{array}}
\makeatletter
\@addtoreset{equation}{section}
\makeatother

\def\CV{{\mathbb V}}

\def\Mark{{\boldsymbol \tau}}
\def\B{{\Omega}} 
\def\Hc{{\mathcal H}}
\def\Te{\mathcal T} 
\def\Proj{ \mathbb S}
\def\Projm{\mathbb S}
\def\M{{\mathcal M}}
\def\Qm{{\widehat{{\mathcal Q}}}}
\def\L{{\mathcal L}}

\def\Mc{{\cal M}}
\def\tr{{\rm tr}\, }
\def\phi{\varphi}

\def\z{{\zeta}}
\def\C{{\mathbb C}}
\def\CP1{{\mathbb C\mathbb P}^1}

\def\Z{{\mathbb  Z}}
\def\R{{\mathbb R}}
\def\a{\alpha}
\def\g{\gamma}
\def\b{\beta}

\def\l{\lambda}

\def\s{\sigma}

\def\o{\omega}
\def\O{\Omega}
\def\p{\partial}

\def\Ch{{\widehat{\mathcal C}}}
\def\Bh{{\widehat{B}}}

\def\proj{{S}}
\def\Q{{\mathcal Q}}

\def\spin{\kappa}

\def\pa{\partial}

\def\gt{{\tilde{\gamma}}}

\def\vh{{\hat{v}}}
\def\wh{\widehat{w}}

\def\at{\tilde{a}}
\def\bt{\tilde{b}}

\def\dim{{\rm dim}}
\def\fh{{\hat{f}}}
\def\wh{\widehat}

\def\bh{\hat{b}}
\def\Bh{{\widehat{B}}}

\def\gh{\widehat{g}}

\def\bh{\widehat{b}}

\def\Mgn{\mathcal{Q}_{g}}

\def \Im {\operatorname{Im}}
\def\f{\frac}
\def\la{\label}
\def\Mod{{\mathcal M}}
\def\Tc{{\mathcal T}}

\def\Ld{{\Lambda}}  
\def\H{H}  
\def\h{h}  
\def\bc{q} 
\def\qd{Q} 
\def\Scal{{\mathcal S}} 
\def\Ldm{{\Ld}_-}
\def\Ldp{{\Ld}_+}
\def\Ldt{{{\Ld}_3}}
\def\VF{{\mathbf V}} 

\def\ocan{\omega _{{can}}}
\def\tcan{\theta_{{can}}}

\def\res{{\rm res}}
\def\det{{\rm det}}
\def\ov{\overline}

\def \1{{\mathrm I}}

\def \bd{\begin{definition}}
\def\ed{\end{definition}}

 \renewcommand{\familydefault}{\sfdefault}
\begin{document}

%
%
%
%

\baselineskip 14pt plus 1pt minus 1pt

\vspace{0.2cm}
\begin{center}
\begin{Large}
\fontfamily{cmss}
\fontsize{17pt}{27pt}
\selectfont
\textbf{
Symplectic geometry of the moduli space of projective structures in homological coordinates}
\end{Large}\\
\bigskip
M. Bertola$^{\dagger\ddagger\diamondsuit}$\footnote{Marco.Bertola@concordia.ca, mbertola@sissa.it},  
D. Korotkin$^{\dagger\ddagger}$ \footnote{Dmitry.Korotkin@concordia.ca},
C. Norton$^{\dagger\ddagger}$ \footnote{Chaya.Norton@concordia.ca}
\\
\bigskip
\begin{small}
$^{\dagger}$ {\it   Department of Mathematics and
Statistics, Concordia University\\ 1455 de Maisonneuve W., Montr\'eal, Qu\'ebec,
Canada H3G 1M8} \\
\smallskip
$^{\ddagger}$ {\it  Centre de recherches math\'ematiques,
Universit\'e de Montr\'eal\\ C.~P.~6128, succ. centre ville, Montr\'eal,
Qu\'ebec, Canada H3C 3J7} \\
\smallskip
$^{\diamondsuit}$ {\it  SISSA/ISAS,  Area of Mathematics\\ via Bonomea 265, 34136 Trieste, Italy }\\
\end{small}
\vspace{0.5cm}
{\bf Abstract} \end{center}
{
 We introduce a natural symplectic structure on the moduli space of quadratic differentials with simple zeros and describe its Darboux coordinate systems in terms of  so-called homological coordinates. We then show that this structure coincides with the canonical Poisson structure on the cotangent bundle of the moduli space of Riemann surfaces, and therefore the homological coordinates provide a new system of Darboux coordinates. 
We  define a natural family of commuting "homological flows" on the moduli space of quadratic differentials and find the corresponding action-angle variables. 

The space of projective structures over the moduli space can be identified with the cotangent bundle upon selection of a reference projective connection that varies holomorphically and thus can be naturally endowed  with a symplectic structure. Different choices of  projective connections of this kind (Bergman, Schottky, Wirtinger) give rise to equivalent symplectic structures on the space of projective connections but different symplectic polarizations: the corresponding generating functions are found. 
We also study the monodromy representation of the Schwarzian equation associated with a projective connection, and we show that the natural symplectic structure on the the space of projective connections induces the Goldman Poisson structure on the character variety. Combined with results of Kawai, this result  shows the symplectic equivalence between the embeddings of the cotangent bundle into the space of projective structures given by the  Bers and Bergman projective connections.
}

\vspace{0.7cm}

\vskip 15pt

\tableofcontents
%
%
%
\section{Introduction and results}
 The monodromy map of the Schwarzian equation is a classical topic which has attracted 
the attention of many 
authors. For a comprehensive review see \cite{GaKaMa} and classical works \cite{Tyurin1978,Earle,Hejhal,Gunning}.

The study of its symplectic aspects was initiated by  S.Kawai \cite{Kawai} who established a relationship 
between the canonical symplectic structure on the cotangent bundle of the Teichm\"uller space
and    Goldman's bracket \cite{Gold86} for the traces of  monodromy matrices. Kawai's results were  recently used in  the context of Liouville theory (see \cite{NeRoSa,Teschner2012} for details).

The principal goal of this paper is to present an alternative and explicit approach to the symplectic geometry of this monodromy map.
To describe our results we introduce
 the following bundles of complex dimension $6g-6$
over the moduli space $\M_g$ of Riemann surfaces of genus $g$:
\vspace{-5pt}
\begin{packed_enum}
\item[--] The cotangent bundle $T^*\Mc_g$, which can be identified with the  space $\Q_g$
 of all holomorphic quadratic differentials. A point of $\Q_g$ is a pair $(\CC,\qd)$ with $\CC$ a closed Riemann surface of genus $g$ and $\qd$ a holomorphic quadratic differential on $\CC$. The universal cover $\Qm_g$ of $\Q_g$ is the cotangent bundle $T^*\Te_g$ over the Teichm\"uller space $\Te_g$. 
 \item[--] The open subset  $\Q^0_g \subset \Q_g$ of holomorphic quadratic differentials with simple zeros  over non-hyperelliptic Riemann surfaces of genus $g$ (we exclude the hyperelliptic locus to avoid unnecessary technical details).
 \item[--] The moduli space of projective connections  $\Proj_g$, which is an affine bundle over $\M_g$. The fiber of  $\Proj_g$ is the affine vector space of holomorphic  projective connections of fiberwise dimension $3g-3$.
\item[--] The moduli space of  projective structures. This moduli space can be identified with the space of conjugacy classes of 
 ``non-elementary" (see \cite{GaKaMa})  $PSL(2,\C)$ representations of the fundamental group of $\CC$. The monodromy map of the Schwarzian equation
\be
{\cal S}(f,\cdot)= S(\cdot)
\la{Schwint0}
\ee
where $\mathcal S$ is the Schwarzian derivative and $S$ is a given projective connection, 
assigns a  non-elementary  $PSL(2,\C)$ representation to the pair (Riemann surface 
$\CC$, projective connection $S$ on $\CC$).
This $PSL(2,\C)$ representation can be  ``lifted'' to an $SL(2,\C)$ representation \cite{GaKaMa}. The space of conjugacy classes of such representations will be denoted $\CV_g$  and is called the {\it character variety}.
\end{packed_enum}

For any point $(\CC,\qd)\in \Q_g^0$, the canonical 2-sheeted covering curve $\Ch$  is defined by the locus of the equation $v^2=\qd$ in $T^*\CC$. This cover is branched at  the $4g-4$ zeros of $\qd$, and its genus is  $\wh g=4g-3$. Furthermore $v$ is a holomorphic differential on $\Ch$ with zeros of order $2$ at the branch points. 
Denote by $\mu$ the natural holomorphic involution on $\Ch$. The differential   $v$ is anti-symmetric under the action of $\mu$. The first homology of $\Ch$ (over $\R$ for convenience) decomposes into the direct sum $H_+\oplus H_-$ of the even and odd eigenspaces with respect to the action of $\mu$; $\dim H_+=2g$ and $\dim H_-=6g-6$.

 The integrals $\oint_s v$, for   $s$  ranging in a basis of $H_-$, give local holomorphic coordinates in $\Q_g^0$ which are called {\it homological coordinates} \cite{DuaHub,KonZor}.
We define a natural Poisson bracket for these coordinates,  also called  {\it homological};
\be
\le\{
\oint_s v, \oint_{\wt s} v
\ri\} := s\circ \wt s.
\label{poissint0}
\ee
Since the intersection pairing on $H_-$ is nondegenerate, \eqref{poissint0} defines a (homological) symplectic structure; Darboux coordinates are obtained by choosing any symplectic basis $(s_1,\dots s_{3g-3}; s_1^\star, \dots, s_{3g-3}^\star)$ with $s_j \circ s_k^\star = \delta_{jk}$. Denoting by $\P_{s_j}, \P_{s_k^\star}$ the corresponding homological coordinates, the symplectic form corresponding to the Poisson structure (\ref{poissint0}) takes the form:
\be
\o= \sum_{j=1}^{3g-3} \d \P_{s_j}\wedge \d \P_{s_j^\star}\;.
\la{ohom}
\ee

On the other hand, the moduli spaces $\Q_g$ and $\Qm_g$ carry a canonical (complex-analytic) symplectic structure. Given a system of local coordinates $q_1,\dots,q_{3g-3}$ on  
$\Mc_g$ (avoiding orbifold points) or on $\Te_g$, the differentials $\d q_1,\dots,\d q_{3g-3}$ form a basis of the cotangent space, and any holomorphic quadratic differential $\qd$ can be represented as 
$\qd=\sum_{i=1}^{3g-3} p_i \d q_i$. 
The canonical symplectic 2-form on $\Q_g$ and $\Qm_g$ is defined by
\be
\o _{{can}}=\sum_{i=1}^{3g-3} \d p_i\wedge \d q_i\;.
\la{ocanint}
\ee
This form is  independent of the choice of the local coordinates $\{q_i\}$. 
The canonical symplectic potential is defined by $\theta _{{can}}=\sum_{i=1}^{3g-3} p_i \d q_i$ and $\d \theta _{{can}} = \omega _{{can}}$.

Our first result is given by the following theorem (see Thm. \ref{taut}); 
\begin{theorem}\la{th1int}
The forms (\ref{ohom}) and (\ref{ocanint}) coincide, that is 
$$\o=\o _{{can}}\;$$
on the intersection of their respective domains of definition. In terms of the  homological coordinates the canonical symplectic potential satisfies
\bea
\theta _{{can}}=\frac 1 2 \sum_{j=1}^{3g-3} (\P_{s_j} \d \P_{s_j^\star}- \P_{s_j^\star} \d \P_{s_j})\;.
\la{thetacanint}
\eea
\end{theorem}

 Under a change of Torelli marking the homological coordinates $\{\P_j\}$ undergo a complicated nonlinear transformation due to the fact that the zeros  of the quadratic differential 
 $Q^\sigma=Q+\f{1}{2}(S_B-S_B^\s)$ have different locations,  and thus the equation $v^2=\Q^\sigma$ defines another canonical double covering curve   $\Ch^\sigma$. 
Theorem \ref{th1int} implies that
\be
\sum_{j=1}^{3g-3} \d \P_{s_j}\wedge \d \P_{s_j^\star}=
\sum_{j=1}^{3g-3} \d \P^\sigma_{s_j}\wedge \d  \P^\sigma_{s_j^\star}
\la{strstr}\ee
where $\{\P_{s_j}, \P_{s_j^\star}\}$ and $\{ \P^\sigma_{s_j}, \P^\sigma_{s_j^\star}\}$  are  homological coordinates corresponding to two symplectic bases in $H_-(\Ch)$ and $H_-(\Ch^\sigma)$, respectively. A direct verification of the identity (\ref{strstr}) would be technically non-trivial.
\bigskip

The Prym matrix $\Pi$  (see (\ref{Prymmat})) of the canonical cover $\Ch$  plays the main role in the change from canonical to homological coordinates on $T^*\Mc_g$.
In Sect. \ref{SecPrym} we find a  Lie--commuting basis of vector fields  $\{\VF_j\}_{j=1}^{3g-3}$  spanning  the ``vertical'' foliation of $T^*\M_g$. These vector fields are expressed in terms of homological coordinates, and the Prym matrix is essential for their construction. 
 Their commutativity  implies in particular a system of PDE's for the Prym matrix with respect to homological coordinates (equation (\ref{sysPi}) below) 
which we were unable to find in the existing literature.

\bigskip
Any  polarization  for the canonical Poisson structure on $\Q_g$ can be used to define a completely integrable system on $\Q^0_g$.  In view of Thm. \ref{th1int},  a 
natural system of commuting Hamiltonians is given by  
\be
   H_j=\f{1}{4\pi} {\P_{s_j}^2} \;,\hskip0.7cm j=1,\dots,3g-3.
   \la{HAint}
\ee

Let $\Q_H$ be the joint level set of the $H_j$'s;
we show  in Sect. \ref{secthomo} that these Hamiltonians are, in fact, action variables, namely, $H_j = \oint_{\varpi_j} \theta _{{can}}$ with $\varpi_j$ a non-contractible loop in $\Q_H$. 
 Similarly to Goldman's flows \cite{Gold86} these "homological flows" commute amongst themselves and are holomorphic in the moduli.
\bigskip

We now discuss the relationship with the space $\Proj_g$ of projective connections. 
Let $S_0$ be a holomorphic projective connection on $\CC$ that depends holomorphically on the moduli in a neighborhood of $\CC\in \M_g$; we can identify the fibers of $\Proj_g$ and $\Q_g$ over this neighborhood 
by the map 
\be
\qd\to S=S_0+ 2\qd.
\la{qtoS}
\ee
In general, a global definition  of $S_0$  requires a suitable marking of $\CC$:  Teichm\"uller, Torelli or Schottky.
The map (\ref{qtoS})  induces a symplectic structure on $\Proj_g$ from the canonical symplectic structure  on $\Q_g$  depending on the choice of 
reference projective connection $S_0$ and therefore defines a Lagrangian embedding of the base into $\Proj_g$.
We call two reference projective connections $S_0$ and $S_1$ {\it equivalent} if the corresponding induced
 symplectic structures on $\Proj_g$ coincide. Within the same equivalence class,  different choices of the reference connection 
 define  different Lagrangian embeddings.
 In Sect. \ref{SecCan} we show that 
$S_0$ and $S_1$ are equivalent if and only if there is a holomorphic function $G$ of the moduli (and depending on the appropriate marking), such that   $\d G$ 
corresponds to  the  quadratic differential $ \f{1}{2}(S_0-S_1)$ under the standard identification of $\Q_g$ and $T^*\Mc_g$ (in the rest of the paper we use 
this identification without further comments).
 The function $G$ is  shown to be the generating function of the change of 
 Lagrangian embedding  induced by the change of reference projective connection from $S_0$ to $S_1$,
 i.e. the  function generating the change of corresponding systems of Darboux coordinates.

In this paper we primarily  use  the   Bergman projective connection $S_B$, which depends on a Torelli marking of $\CC$  and transforms under the action of  a symplectic transformation 
\be
\sigma=    \left(\ba{cc} D & C\\B  & A  \ea\right)\;\;\in\;\;{\rm Sp}(2g,\Z)\;.
\label{symint}
\ee
as explained in  Sect. \ref{projconnect} (see also  \cite{Fay73}).

Other projective connections we consider (see Sect. \ref{projconnect} for definitions and details): 
\vspace{-5pt}
\begin{packed_enum}
\item[--] the Schottky 
projective connection $S_{Sch}(\cdot)={\cal S}(w,\cdot)$ where $w$ is the Schottky uniformization coordinate on $\CC$;
\item[--]
the Wirtinger projective connection $S_W$  which is independent of the Torelli marking and is defined on $\M_g$ outside of the divisor where an even spin structure admits holomorphic sections;
\item[--] the Bers projective connection $S_{Bers}$ which was used as the reference connection in \cite{Kawai}. 
\end{packed_enum}
Our next result states the equivalence of Bergman, Schottky and Wirtinger projective connections (Section \ref{SecCan} and in particular Prop. \ref{propGstt}, \ref{propWirt}
 and \ref{propSchot}).
\begin{theorem} \la{th2int}
The Bergman projective connections for different Torelli markings, the Schottky projective connection and the Wirtinger projective connection are equivalent. 
\newline
The generating functions corresponding to the change of Lagrangian embeddings are:
\begin{packed_enum}
\item[--]  For  two Torelli markings  related by  a symplectic transformation \eqref{symint}
$
G=6\pi i \log  {\rm det} (C\Omega+D)\;.
$
\item[--] For the change from Bergman to Schottky projective connections $G=-6\pi i \log F$, where $F$ is the Bowen-Zograf F-function on Schottky space \cite{Bowen,Zograf1990,MacTak}.
\item[--] For the change  from Bergman to Wirtinger projective connections
 $G=-(24\pi i/(2^g+4^g)) \log \Theta$ where $\Theta$ is the product of $2^{g-1}+2^{2g-1}$ theta-constants.
 \end{packed_enum}
\end{theorem}

The  monodromy map of the Schwarzian
equation (\ref{Schwint0}) defines a $PSL(2,\C)$ representation of $\pi_1(\CC)$. 
In Sect. \ref{secmon} we define a suitable $SL(2,\C)$ lift  (see also \cite{GaKaMa}), and hence a map to $\CV_g$. 
Under this monodromy map, the symplectic structure on $\Proj_g$ induces a nondegenerate Poisson structure on  $\CV_g$.

The developing map of a projective connection can  be written as the ratio of two solutions of  a scalar second order ODE on $\Ch$ (Sect. \ref{Seceq}) 
\be
\psi_{zz}  -u(z) \psi=0\;,
\label{matrix1in}
\ee
where  $z(x)= \int_{x_1}^x v$ (here $x_1$ is a selected zero of $\qd$) is a local coordinate on $\CC$.  The  coefficient  $u$ in \eqref{matrix1in} is given by
\be
u=-\f{S_B-S_v}{2\qd}-1
\la{uint}
\ee
where $S_v(\cdot)=\Scal\le(\int_{x_1}^x v, \cdot\ri)$, and it is a meromorphic projective connection on $\CC$. 

We will show (Thm. \ref{thmPoissonu}) that  the Poisson bracket of
$u(z)$ and $u(\zeta)$ is given by
\be
\f{4\pi i}{3}\{u(z), u(\zeta)\}=\L_z \left[ h^{(\zeta)}(z)\right]- \L_\zeta \left[h^{(z)}(\zeta)\right]\;,
\label{Poissonuin}
\ee
where the differential operator 
$$\L_z:= \frac 1 2 \pa_z^3 - 2u(z) \pa_z  -u_z(z)$$ is known as the  Lenard operator in the theory of KdV equation
 \cite{BaBeTa}; $h^{(y)}(x):=\int_{x_1}^x B^2(x,y)/(\qd(x)\qd(y))$ is a meromorphic function on $\CC\times\CC$ and $B$ is the fundamental normalized bidifferential on $\CC\times\CC$ with a
 second order pole on the diagonal (Sect. \ref{projconnect}).  The Poisson bracket  (\ref{Poissonuin}) for $u(z)$  implies the following Theorem  (Thm. \ref{sympGold}, Sect. \ref{SecPoi}):

\begin{theorem}\la{th3int}
Given $\g, \wt \g\in \pi_1(\CC, x_0)$, the Poisson bracket between traces of the corresponding monodromy matrices $M_\g, M_\gt$  of the equation (\ref{matrix1in}) is given by:
\be
\{\tr M_\g,\;\tr M_\gt\}= \f{1}{2}\sum_{p\in \g\cap \gt} \le(\tr M_{\g_p\gt}-\tr M_{\g_p\gt^{-1}}\ri)
\label{Goldmanint}
\ee  
where $\g_p\gt$ and $\g_p\gt^{-1}$ are two ways to resolve the intersection point $p$ to get two new contours  $\g_p\gt$ and $\g_p\gt^{-1}$ for each $p\in \g\cap \gt$. 
\end{theorem}
The right side of \eqref{Goldmanint} is the Goldman Poisson bracket on the (complexified) $SL(2,\R)$  character variety \cite{Gold86}, and hence the homological and Goldman Poisson structures coincide.

The same bracket on $\CV_g$ was deduced by S.Kawai \cite{Kawai} in 1996 using a Bers projective connection $S_{Bers}$ as a reference connection; the
coincidence of the brackets on $\CV_g$ shows the equivalence of the Bergman and Bers projective connections, which would be  elusive to establish directly.
Therefore, Theorem \ref{th3int}, together with \cite{Kawai} implies the following 
\begin{corollary}
There exists a generating function $G$  on $\Te_g$ such that 
\be
d G=\f{1}{2} (S_B-S_{Bers}^{\CC_0,\eta})\;.
\la{dG}\ee
\end{corollary}
The  function $G$ satisfying (\ref{dG}) is conjecturally given by the derivative of a quasi-Fuchsian generalization of Selberg's zeta-function at 1 (see \eqref{conjfor}).
\bigskip 

\noindent{\bf Acknowledgments.} We thank Andrew McIntyre for several clarifications.
DK thanks Tom Bridgeland and Peter Zograf for interesting discussions. Research support  of MB and DK was provided by NSERC and  FQRNT.

%
%
\section{Projective connections on Riemann surfaces}
\la{projconnect}

We recall that a {\it  holomorphic projective connection} $S$  on a Riemann surface $\CC$ is a collection of representative holomorphic functions in each local coordinate chart, and its representative functions $S(\z), S(\xi)$ in different overlapping charts with coordinates $\z,\xi$ are related by the affine transformation rule
\be
\label{1312}
S(\z)  = \le(\frac {\d \z}{\d \xi}\ri)^2  S(\xi)  + \mathcal S(\z,\xi)\ .
\ee
Here  $\Scal$ is the Schwarzian derivative
\be
 \Scal(\z,\xi)=\le(\frac{\z''}{\z'}\ri)'-\frac 1 2\le(\frac{\z ''}{\z '}\ri)^2\ ,
\la{Schwint}
\ee
the prime denoting the derivative with respect to $\xi$.

\paragraph{Torelli markings and Bergman projective connections.}
Let $\CC$ denote a Riemann surface of genus $g$ and consider the fundamental group $\pi_1(\CC, x_0)$ based at a point $x_0\in \CC$.
We introduce a standard
 set of generators  $\{\a_i,\b_i\}_{i=1}^g$
of   $\pi_1(\CC,x_0)$ satisfying the relation
\be
\prod_{i=1}^g \a_i\b_i\a_i^{-1}\b_i^{-1}=id\;.
\label{relation}
\ee

The corresponding cycles in the homology group $H_1(\CC,\Z)$  will be denoted by $a_i, b_i$ respectively, and the set $\{a_i,b_i\}_{i=1}^g$ is a canonical symplectic basis of cycles in $H_1(\CC,\Z)$, called a {\it Torelli marking} of $\CC$. 
We denote by ${\bf v}=(v_1,\dots,v_g)$ the dual basis of normalized holomorphic differentials satisfying
\be 
\oint_{a_i} v_j=\delta_{ij},\qquad \Omega_{ij}=\oint_{b_i} v_j\;\;\;\hskip0.7cm i,j=1,\dots,g
\label{bperC}
\ee
where the matrix $\Omega$ is called the {\it period matrix} of $\CC$.

Depending on the  Torelli marking $\Mark$ we introduce the {\it fundamental meromorphic bidifferential} \cite{Fay73} $B (x,y)$ on $\CC$ defined by the following properties;
\begin{itemize}
\item  $B(x,y)$ is holomorphic on $\CC\times \CC$ except for a pole of order two along the diagonal with biresidue 1.
\item it is symmetric $B(x,y)= B(y,x)$; 
\item its $a$--periods vanish,  $\oint_{y\in a_i} B(x,y) = 0$, $\forall i=1,\dots g$;

\end{itemize}

If the points $x,y$ belong to the same chart of the local coordinate $\zeta$, then $B(x,y)$ has the following expansion near the diagonal
\be
B(x,y)=\f{d\zeta(x) d\zeta(y)}{(\zeta(x)-\zeta(y))^2} +\f{1}{6}S_B\left(\f{\zeta(x)+\zeta(y)}{2}\right)+ \mathcal O((\zeta(x)-\zeta(y))^2).
\label{Bdiag}
\ee
The term $S_B$ in \eqref{Bdiag} (also depending on Torelli marking of $\CC$) transforms as  a  projective connection  under a change of the local coordinate; it is called the "Bergman projective connection". The Bergman projective connection is defined at any point of the Torelli space.
The bidifferential $B$ can be written explicitly in terms of theta functions (see \cite{Fay73}, p. 20):
$$
B(x,y) = \d_x \d_y \ln \theta\le(\int_x^y {\bf v} - {\bf f}\ri)
$$
where $\theta$ is the Riemann theta function and ${\bf f}\in \C^g$ is any point in the smooth locus of the theta-divisor 
(the result is independent of this choice). This allows one to express $S_B$  in terms of  theta-functions (\cite{Fay73}, p.19), but the explicit expression will not be needed in this paper

\paragraph{Wirtinger  projective connection.}
\label{Wirtproj}
The Wirtinger projective connection (see \cite{Tyurin1978} and \cite{Fay73}, p. 22)  is defined by:
\be
\label{Swirt}
S_W=S_B+ \f{48\time 4\pi i }{2^g+4^g}\sum_{i,j=1}^g v_i v_j\f{\p}{\p\O_{ij}}\log \left(\prod_{\beta\;\; even}\theta[\beta](0)\right)
\ee
where the product is taken over all even theta-characteristics. The Wirtinger projective connection is independent of the Torelli marking;  it is defined at any point of the moduli space 
$\Mod_g$ except along the divisor where one or more of the  theta-constants in \eqref{Swirt} vanish.

\paragraph {Schottky projective connection.}
Another projective connection that depends  holomorphically  on moduli is the Schottky projective connection. It is defined as the Schwarzian derivative
\be
\label{Schottproj}
S_{Sch}(\cdot)={\mathcal S}(w,\cdot)
\ee
where $w$ is the Schottky uniformization coordinate. 
The definition of the Schottky uniformization requires fixing $g$ generators $\a_1,\dots,\a_g$ of the fundamental group (see for example \cite{MacTak} for details). The moduli space of these "Schottky-marked" Riemann surfaces is called the Schottky space. To each point of the Schottky space one can assign a Schottky
group $\Gamma_S$ (a subgroup of $PSL(2,\C)$) such that the Riemann surface $\CC$ is identified with the quotient $\C/\Gamma_S$, and $w$ is then the uniformizing coordinate in a fundamental domain of $\C/\Gamma_S$. The Schottky group is the monodromy group of the Schwarzian equation \eqref{Schottproj}.

\paragraph {Bers projective connection. }
\label{Bersproj}
 Fix a reference point  $(\CC_0, \{\a_i,\b_i\})$ in $\Te_g$, and let $\Gamma_0\subset SL_2(\R)$ be the Fuchsian group  corresponding to $\CC_0$. Recall that the Fuchsian projective connection on $\CC_0$ is given by the Schwarzian derivative ${\mathcal S}(t,\cdot)$
where $t$ is the Fuchsian  uniformization coordinate of $\CC_0$. 

Choose a holomorphic quadratic differential $\eta(t)$ on $\CC_0$ which can be viewed as an automorphic  function of weight 4 on the upper half-plane with respect to $\Gamma_0$. By solving the Beltrami equation $
\pa_{\ov t} W(t,\ov t)  = \chi_{\mathbb H_+}\Im(t)^4\,\ov {\eta(t) }\pa_t W(t,\ov t)$ ($\chi_{\mathbb H_+}$ is the indicator function of the upper half-plane 
${\mathbb H_+}$), one obtains a {\it quasi-conformal} transformation $W$ which is holomorphic in the lower half-plane $\mathbb H_-$ mapping the $t$-plane to the $t_\eta$-plane.
If  $\sup_{t\in \mathbb H_+} |\Im(t)^4 \eta(t) |< 1$\footnote{The sup is actually a maximum because $|\Im(t)^4 \eta(t) |  $ is  invariant under the action of $\Gamma_0$ and hence it suffices to take the sup on a fundamental domain.}, the map is a diffeomorphism of the Riemann sphere which can be uniquely determined by requiring $W(0)=0$, $W(1)=1$ and $W(\infty)=\infty$.  The automorphism group $\Gamma_{\CC_0,\eta}:= W \Gamma_0 W^{-1}$ of the $t_\eta$-plane is a group of M\"obius transformations of the quasi-Fuchsian type;  $W$ maps $\R$ into  a simple Jordan curve, and $\Gamma_{\CC_0,\eta}$ 
acts properly discontinuously both in the interior and in the exterior of this curve. Therefore $Int(W(\R))/\Gamma_{\CC_0,\eta}$ defines a new Riemann surface $\CC_\eta$, whereas $Ext(W(\R))/\Gamma_{\CC_0,\eta}$ is equivalent to a reflection of $\CC_0$ (this is "Bers' simultaneous uniformization" \cite{Bers59}).

By construction the coordinate $t_\eta$ is the developing map with monodromy group $ \Gamma_{\CC_0,\eta}$. The Schwarzian 
derivative ${\mathcal S}(t_\eta,\cdot)$ defines a holomorphic projective connection on $\CC_\eta$ depending on the reference Riemann surface $\CC_0$ (as well as on the marking). 
This projective connection depends holomorphically on the moduli of  the quadratic differential $\eta$, and non-holomorphically on the moduli of $\CC_0$ . We shall denote such a projective connection by $S_{Bers}^{\CC_0,\eta}$. The 
group $\Gamma_{\CC_0,\eta}$ is the monodromy group of the 
Schwarzian equation ${\mathcal S}(t_\eta,\cdot)=S_{Bers}^{\CC_0,\eta}(\cdot)$.

\medskip

\paragraph{Transformation of the canonical bidifferential under the change of Torelli marking.}
Both   $B(x,y)$ and $S_B$ depend on a Torelli marking of  $\CC$ . Given a matrix $\sigma \in Sp(2g,\Z)$ (\ref{symint}),  a new  basis 
in $H_1(\CC,\Z)$ is defined by
\be
 \left(\ba{c} a_i^\s\\ b_i^\s \ea\right)=\s
\left(\ba{c} a_i\\ b_i \ea\right)\;.
\label{symptrint}
\ee
Then the period matrix \eqref{bperC}, the fundamental bidifferential and the Bergman projective connection transform as follows  (\cite{Fay73}, p.21);
\be
\Omega^{\s}  = (A\Omega+ B)(C\Omega+D)^{-1}\;,
\label{transOC}
\ee
\be
B^{\s} (x,y)  = B(x,y) 
-\pi i \sum_{ j\leq k\leq g} \le(v_j(x) v_k(y) + v_k(x) v_j(y) \ri) \frac {\pa \ln \det (C\O+D)}{\pa \O_{jk}} \;,
\label{transB}
\ee
\be
S_{B}^{\s}  = S_B -12\pi i 
{\sum_{ j\leq k\leq g} v_jv_k \frac \pa {\pa \O_{jk}} \ln \det (C \O+D)
}\;,
\label{transSB}
\ee
where $(v_1,\dots,v_g)$  are normalized  as in \eqref{bperC}. 
\section {The canonical double cover}

\label{Cancov}

Let $\Q_g^0$ denote the open subset of 
the space $\Q_g$ consisting of quadratic differentials
with simple zeros on non-hyperelliptic Riemann surfaces of genus $g$. The dimension of $\Q_g^0$ equals the dimension of $\Q_g$ which is $6g-6$. 
To each point $(\CC,\qd)\in \Q_g^0$ we associate a new Riemann surface $\Ch$ (called the "canonical cover") defined by the equation  
\be
v^2=\qd
\label{defv}
\ee
in  the cotangent bundle $T^\star \CC$. 
The canonical projection $\pi: \Ch \mapsto \CC$ has branch points at the zeros $\{x_1,\dots, x_{4g-4}\}$ of $\qd$. We denote by $\mu :\Ch \to \Ch$ the holomorphic involution on $\Ch$ which interchanges the sheets. 

\subsection{Geometry of the double cover}

The genus of $\Ch$ equals $\hat{g}=4g-3$, and the preimages of the zeros under the canonical projection are {\it double} zeros of the differential $v$ on $\Ch$: indeed, in a suitable coordinate $\zeta$ we have $\qd = \zeta \d \zeta^2$, and the local coordinate on $\Ch$ is $\xi = \sqrt{\zeta}$ which implies $v = \xi \d (\xi^2) = 2 \xi^2 \d \xi$. 

\begin{figure}[ht]
\centering
\resizebox{0.3\textwidth}{!}{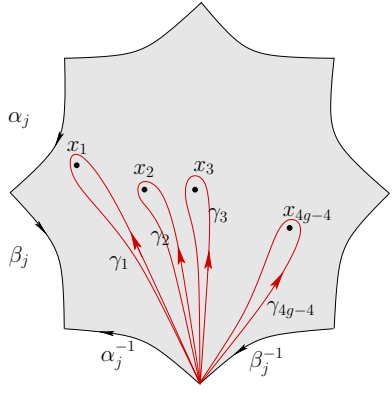}\\
\captionof{figure}{Fundamental polygon $\CC_0$ with loops $\g_i$ around zeros $x_i$}
\label{genpi1}
\end{figure}

Fix a point $x_0\in \CC$ and consider the fundamental group $\pi_1(\CC\setminus\{x_i\}_{i=1}^{4g-4}, x_0)$. We choose a standard 
 set of generators
 $\{\a_i,\b_i,\g_j\}$, $i=1,\dots,g$,  $j=1,\dots 4g-4$ where the contours   $\g_j$ are freely homotopic to small circles around the zeros  $x_{1},\dots x_{4g-4}$ of $\qd$. The generators are chosen so that they satisfy the relation:
\be
\g_{4g-4}\cdots\g_{1}\prod_{i=1}^g \a_i\b_i\a_i^{-1}\b_i^{-1}=id\;.
\label{relpunct}
\ee
The holonomy  of $v$ on $\CC\setminus\{x_i\}_{i=1}^{4g-4}$  defines a homomorphism 
\be
\mathfrak h: \pi_1(\CC\setminus\{x_i\}_{i=1}^{4g-4}, x_0) \mapsto \Z_2\;,
\ee
where   necessarily  $\mathfrak h(\g_j)=-1$,  while $\mathfrak h(\a_i)$ and $\mathfrak h(\b_i)$ depend on the choice 
of generators. 
The following simple lemma guarantees that for some choice of $\{\a_i,\b_i\}$, the homomorphism $\mathfrak h$ is the identity on $\a_i$ and $\b_i$ for $i=1,\dots,g$.
\begin{lemma}
\label{S2lem}
The generators $\{a_i,\b_i,\g_j\}$ of $\pi_1(\CC\setminus\{x_i\}_{i=1}^{4g-4},x_0)$ satisfying (\ref{relpunct}) can be chosen such that 
\be
\mathfrak h(\a_i)= \mathfrak h(\b_i)=
1\;\;,\hskip0.7cm \mathfrak h(\g_j)= 
{-1}\;.
\label{geninS2}
\ee
\end{lemma}
\noindent {\bf Proof.} 
The relation (\ref{relpunct}) implies $\mathfrak h(\g_{4g-4}\dots \g_1)= 
{1}$ because the group $
{\Z_2}$ is Abelian. 
Dissect $\CC$ along a system of contours corresponding to the generators $\a_i,\b_i$;
the result of this dissection is  a
 fundamental polygon $\CC_0$ with the zeros of $\qd$ lying 
 in its interior and
the generators $\g_i$ ordered as shown in Fig. \ref{genpi1}. 
The points $x_i$ are branch points, and therefore the sheets of the cover $\Ch$ are always interchanged along $\g_i$, i.e. $\mathfrak h(\g_i)=-1$ for each $\g_i$. Notice that  the number of branch points is even, and this is compatible with the relation $\prod_i \mathfrak h(\g_i)=1$.

\begin{figure}[ht]
\centering
\resizebox{0.3\textwidth}{!}{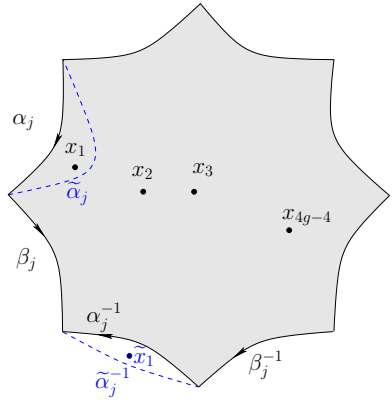}
\\
\captionof{figure}{Transformation of generator $\a_j$ to $\tilde{\a}_j$ such that $\mathfrak h(\tilde{\a}_j)= 
{-}\mathfrak h(\a_j)$}
\label{genpi2}
\end{figure}

Suppose that for some $j$ we have $\mathfrak h(\a_j)= 
-1$; the generator $\a_j$ can then be deformed within its own homotopy class in $\pi_1(\CC, x_0)$  by "moving" $\a_i$ through any zero, for example through $x_1$, see Fig. \ref{genpi2} while changing the other generators accordingly so as to preserve (\ref{relpunct}). This elementary transformation does not modify the images of the other generators under $\mathfrak h$ because the group $\Z_2$ is Abelian. 

The result is a new 
path $\tilde{\a}_j$ satisfying  $\mathfrak h(\tilde{\a}_j)= 1$. 
We then repeat this procedure as necessary until all the generators satisfy (\ref{geninS2}). \QED

\smallskip
\medskip

\begin{figure}[ht]
\centering
\resizebox{0.4\textwidth}{!}{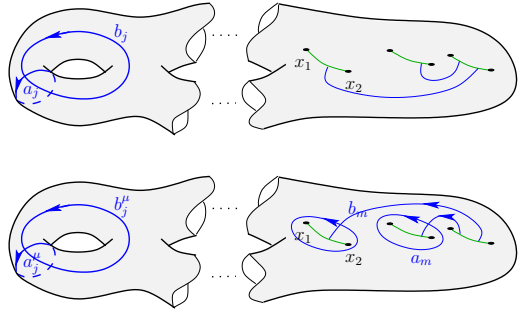}
\\
\captionof{figure}{ Choice of canonical basis of cycles  on the canonical cover $\Ch$}
\label{canbasCh}
\end{figure}

\begin{figure}[ht]
\centering
\resizebox{0.4\textwidth}{!}{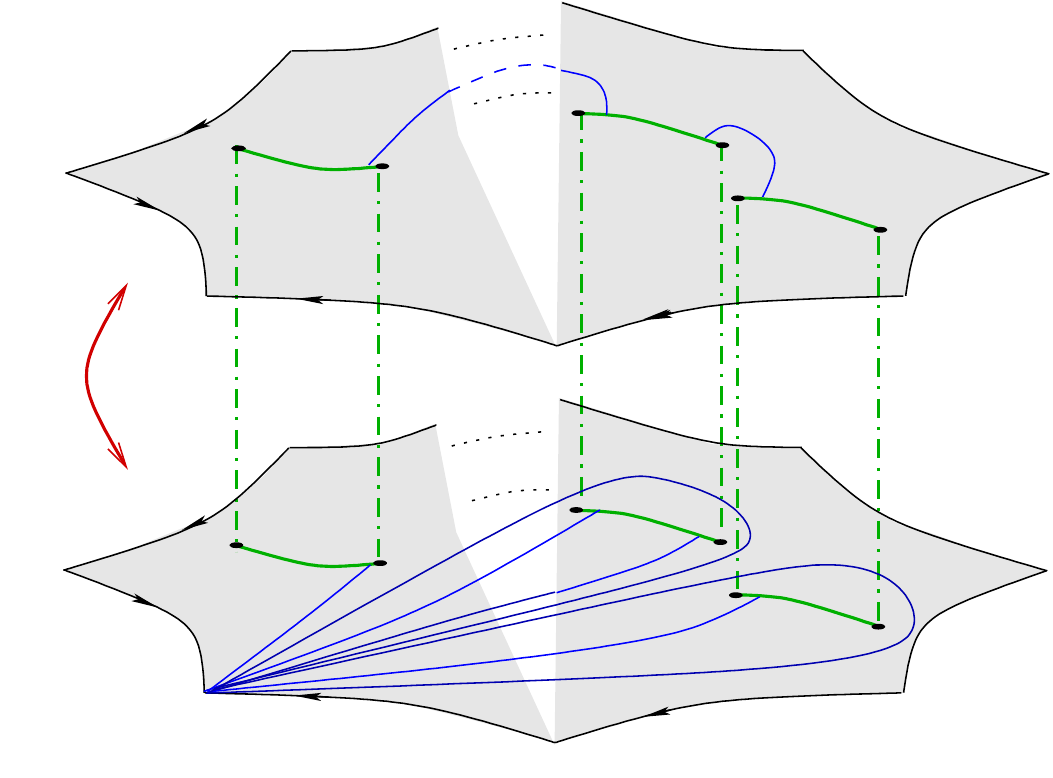}
\captionof{figure}{Two-sheeted cover $\tilde{\CC}_0$ of the fundamental polygon $\CC_0$ of $\CC$ with branch points at zeros of $\qd$.}
\label{coverfund}
\end{figure}

Let $\CC_0$ denote the fundamental polygon obtained by dissecting $\CC$ along the generators $\a_j , \b_j$ satisfying \eqref{geninS2}. Choose a 
basepoint $x_0\not \in \{x_1,\dots, x_{4g-4}\}$ for this dissection. Then $\pi^{-1} (\CC_0)$ can be represented (see Fig. \ref{coverfund}) as the two-sheeted cover $\wt\CC_0$ of $\CC_0$ with branch points $x_1,\dots,x_{4g-4}$. The assumption that the ${\Z_2}$ holonomy
around the  loops $\{\a_i,\b_i\}$ is  trivial implies that we can choose  $2g-2$ branch cuts between $x_i$ which do not intersect the boundary of $\CC_0$.
The fundamental polygon $\Ch_0$ of $\Ch$ is obtained by additionally cutting $\wt\CC_0$ along a system of loops
$\{\tilde{\a}_i,\,\tilde{\b}_i\}$ for $i=1,\dots, 2g-3$ shown in Fig.\ref{coverfund} in blue. These loops generate the fundamental group 
of $\wt\CC_0$.  
\paragraph{Even and odd subspaces in the homology and cohomology of $\Ch$.}
The holomorphic  involution $\mu:\Ch \to \Ch$ induces a map with eigenvalues $\pm 1$ on the cohomology of $\Ch$; thus the space $H^{(1,0)}(\Ch)$
 of holomorphic differentials on $\Ch$ can be decomposed into the two eigenspaces $H^\pm$
\be
H^{(1,0)}(\Ch) = H^+(\Ch) \oplus H^-(\Ch) \ ,\qquad \dim \, H^+(\Ch) = g\ ,\qquad \dim \, H^{-}(\Ch) =3g-3\;.
\ee 
The subspace $H^-(\Ch)$ is the space of {\it Prym differentials}; by construction, the differential $v$ satisfies 
$v(x^\mu)= -v(x)$ i.e.  $v\in H^-$.

In turn, the  homology space $H_1(\Ch,\R)$  decomposes in  $H_+(\Ch,\R)\oplus H_-(\Ch,\R)$ 
where ${\rm dim} \, H_+=2g$, ${\rm dim}\,  H_-=6g-6$. 
A canonical basis of $H_1(\Ch,\Z)$ can be chosen as follows:
\cite{Fay73,KorZog2}, see Fig.\ref{canbasCh}
\be
\{a_j,  a_j^\mu, \at_k,\,b_j,  b_j^\mu\,,\bt_k\}\,,\quad j=1,\dots,g,\;k=1,\dots,2g-3\;,
\label{mainbasis}
\ee
where $(a_i,b_i,a_i^\mu,b_i^\mu)$ is a lift  of the canonical basis of cycles on $\CC$ to $\Ch$,
such that 
\be 
\mu_* a_j = a_j^\mu,\quad\mu_* b_j = b_j^\mu,
\quad\mu_* \at_k + \at_k= \mu_* \bt_k + \bt_k  = 0\;.
\ee
We shall denote by 
\be
\label{basisCH}
 {\{\wh v_j, \wh v_{j^\mu}, \wh w_k\}} 
 \ee
  the corresponding basis of normalized Abelian differentials on $\Ch$. 
The differentials $ v_j^+=\wh v_j+\wh v_j^\mu,\; j=1,\dots,g$,
provide a basis in the space $H^+$; these differentials are naturally identified with the normalized holomorphic differentials $v_j$ on 
$\CC$, and therefore, to simplify the notation, they shall often be written $v_j$ instead of $v_j^+$.
 A basis in $H^-$ is given by the $3g-3$ Prym differentials $v_l^-$, where
\be
 v_l^-=\begin{cases} \wh v_{l}- \wh v_{l}^\mu,\quad l=1,\dots, g\;,\\
\sqrt{2}\,\wh w_{l-g}, \qquad\; l=g+1,\dots, 3g-3\;.\end{cases}
\label{Prym}
\ee
The classes in $H_1(\Ch,\R)$ 
\be
a_j^+ = \f{1}{2}(a_j+ a_j^\mu)\;,\quad
b_j^+ = \f{1}{2}( b_j+ b_j^\mu)\;,\quad j=1,\dots,g,
\label{abp}
\ee
with intersection index 
\be
a_j^+\circ b^+_k=\f{1}{2}\delta_{jk}
\ee
form a  basis in $H_+(\Ch ,\R)$, whereas the  classes
\be
a_l^- =  \f{1}{2}(a_{l}- a_{l}^\mu) ,\hskip0.5cm
b_l^- = \f{1}{2}( b_{l}- b_{l}^\mu),\quad l=1,\dots,g\;,
\la{abm1}
\ee
\be
a^-_{l}=\f{1}{\sqrt{2}}\at_{l-g},\;\;\; b^-_{l}=\f{1}{\sqrt{2}}\bt_{l-g}\;,\hskip0.5cm\, l=g+1,\dots,3g-3
\label{abm}
\ee
form a basis in $H_-(\Ch, \R)$ with intersection index
\be
a_j^-\circ b^-_k=\f{1}{2}\delta_{jk}\;.
\label{intH-}
\ee
Integrating the differentials $v_j^+$ over the cycles $\{a_k^+\}$ gives $\int_{a_k^+} v_j^+=\delta_{jk}$, while their integrals over the cycles
$\{b_k^+\}$ give the period matrix of $\CC$: $\int_{b_k^+} {v}_j^+=\O_{jk}$. Similarly, integrating the Prym differentials (\ref{Prym}) over the cycles $\{a_j^-\}$
(\ref{abm1}), (\ref{abm}) yields the $(3g-3)\times (3g-3)$ unit matrix, while their integrals over the cycles $\{b_j^-\}$ give a $(3g-3)\times (3g-3)$ matrix called the {\it Prym} matrix $\Pi$;
\be
\Pi_{jk}=\int_{b_k^-} {v}_j^-\ ,\ \ 1\leq j,k\leq 3g-3,\  \qquad \Pi = \Pi^t.
\la{Prymmat}
\ee
The period matrix of the double cover $\Ch$ can be expressed in terms of $\Omega$ and $\Pi$; therefore, knowing $\O$ and $\Pi$ one can reproduce both
the Riemann surface $\CC$ and the quadratic differential $\qd$ \footnote{A dimension count  would suggest that the complete information about $\CC$ and $\qd$ should be contained in the Prym matrix alone, but we are  not aware of a rigorous proof of this fact.}.

We also define a homomorphism ${\bf g}$ from $\pi_1(\CC)$ to $H_-(\Ch,\R)$
as follows:
\be
{\bf g}(\a_i)= a_i^-= \frac 1 2(a_i - a_i^\mu) \;\;,\hskip0.5cm {\bf g}(\b_i)= b_i^- = \frac 1 2 (b_i - b_i^\mu)\;,\hskip0.7cm i=1,\dots,g\;.
\label{homom1}
\ee
By abuse of notation we shall write $\int_\g w$ as opposed to $\int_{{\bf g}(\g)} w$ for $\g\in \pi_1(\CC,x_0)$ and any Abelian differential $w$ on $\Ch$.

\paragraph{Functions and bidifferentials associated to the canonical cover.}

Here we introduce   several functions and bidifferentials which will be used throughout the paper.

\begin{itemize}
\item[--] 
The Abelian integral in $\Ch$
\be
z(x)=\int_{x_1}^x v
\label{coordz}
\ee
defines a local coordinate on $\Ch\setminus \{x_1,\dots,x_{4g-4}\}$.
In a neighborhood of a branch point $x_j$, the local coordinate can be chosen to be
$$
\zeta_j(x)=(z(x)-z(x_j))^{1/3}\;.
$$

\item[--]
The  scalar function on $\CC$:
\be
\bc(x)=\f{1}{2}\f{S_B-S_v}{v^2}(x)
\label{Qx}
\ee 
 where $S_v(\cdot)=\Scal\le(z,\cdot\ri)$ is the meromorphic projective connection defined by the differential $v(x)$ via the coordinate $z$ \eqref{coordz}.
The function $\bc$ has second order poles on $\CC$ at the zeros of $\qd$.
The following closely related function $u$  will   be   used extensively below: 
\be
u(z)=-(\bc(z)+1)\;.
\label{uq}
\ee

\item[--]
The meromorphic function $b:\Ch\times\Ch \to \CP1$, defined by 
\be
b(x,y)=\f{B(x,y)}{v(x)v(y)}\;.
\label{bpm}
\ee
This function is anti-symmetric with respect to the involution $\mu$ acting on either of its arguments and has a pole of second order on the diagonal $x=y$. 
Additionally $b(x,y)$ has simple poles at the
branch points $x_1,\dots,x_{\gh-1}$ on $\Ch$ with respect to each of its arguments.
Outside of the branch points $\{x_i\}$ the asymptotics of $b(x,y)$ in the local coordinate $z$ is as follows (see \ref{Bdiag}):
\be
b(x,t)=\f{1}{(z(x)-z(y))^2}+ \f{1}{3} \bc(x) +\f{1}{6} \bc_z(x) (z(y)-z(x))+\dots\;
\qquad y\to x\;.
\label{bas}
\ee

\item[--]
The meromorphic functions  $f_j: \Ch \to \CP1$ defined by
\be
 f_j(x)=\f{v_j(x)}{v(x)}\;, \hskip0.7cm j=1,\dots,g \;.
\label{deffj}
\ee
These functions are antisymmetric under the involution $\mu$ 
and have at most simple poles at the branch points $\{x_i\}$.
\item[--]
The bidifferential $\H(x,y)$ on $\Ch$ which has a 4th order pole and a residue on the diagonal and no other singularities:
\be
\H(x,y)=\f{B^2(x,y)}{v(x)v(y)}=b^2(x,y) v(x)v(y)\;.
\label{defH}
\ee

\item[--]
The meromorphic function $\h$ on  $\CC\times \CC$ defined by:
\be
\h(x,y) = b^2(x,y)=\f{\H(x,y)}{v(x)v(y)}=\f{B^2(x,y)}{\qd(x) \qd(y)}\;.
\label{defh}
\ee
This function has simple poles at the zeros of $\qd$ with respect to both arguments. The pullback to $\Ch\times \Ch$ defines a function (symmetric under $\mu$) which will be denoted by the same symbol for brevity. The asymptotics of $h$ near the diagonal  outside of the branch points $\{x_i\}$ in the local coordinate $z$ can be deduced from (\ref{Bdiag}) 
\be
\h(x,y)= \f{1}{(z(y)-z(x))^4}+ \f{2}{3} \f{\bc(x)}{(z(y)-z(x))^2} +\f{\bc_z(x)}{3 (z(y)-z(x))}+ O(1)
\qquad y\to x\;.
\label{asymph}
\ee

\item[--]
Given any closed non-self-intersecting contour $\gamma$ on $\CC$, define the following 1-form on $\Ch$ by integrating $\H$ 
along the cycle $\frac 1 2 \pi^{-1}(\gamma)$ with respect to one of its arguments:
\be
\H^{(\gamma)}(x) = \f{1}{2}\int_{\pi^{-1}(\gamma)} H(x,\cdot)
\label{defHg}
\ee
The form $\H^{(\gamma)}(x)$ is holomorphic on $\Ch$ and has a jump $\frac {2i\pi }{3} q' v = \frac{2i\pi}{3} \d q$ along the contours  $\pi^{-1}(\g)$.

\item[--]
The anti-derivative of $h(x,y)$ on $\Ch$ with respect to one of its variable with the initial point of integration $x_1$:
\be
\h^{(y)}(x) = \int_{x_1}^x \h(y,\cdot) v(\cdot)\;.
\label{hxy}
\ee
Notice that the differential $\h(y,\cdot) v(\cdot)$ is non-singular at $x_1$. For fixed $x$ the integral $\h^{(y)}(x)$ is a non-singlevalued 1-form 
on $\Ch$ with respect to $y$.  (In particular, it has logarithmic singularities at $y=x$, $y=x^\mu$ and $y=x_1$.)
\end{itemize}
\subsection{Variational formulas}
\label{varformulas}
Let $\{s_i\}_{i=1}^{ 6g-6}$ be a basis in 
 $ H_-(\Ch,\Z)$. The periods $ \P_{s_i} = \oint_{s_i } v$ form a system of local coordinates, called "homological coordinates", on  $\Mgn^0$ \cite{DuaHub,KonZor}.
The intersection pairing on $H_-(\Ch, \R) \times H_-(\Ch, \R)\to \R$ is non-degenerate, and therefore one can define a dual basis 
$\{s_i^*\}_{i=1}^{ 6g-6}$ by the condition 
$$
s_i^* \circ s_j = \delta_{ij}\;.
$$

Variational formulas describe the dependence of the period matrix, holomorphic differentials and the canonical bidifferential on the homological coordinates.  
\begin{proposition}
The variation of $f_j = \frac {v_j}v$ with respect to $\P_{s_i}:= \oint_{s_i } v $, keeping $z(x)$ constant, is given by the following formula:
$$
\frac{\p f_j(x)}{\p\P_{s_i}} \Big|_{z(x)={\rm const}}= \f{1}{4\pi i}\int_{s_i^*} \frac{ B(x,t) v_j(t)}{v(t)}\;,
$$
or, equivalently, 
\be
\frac{\p f_j(x)}{\p\P_{s_i} } \Big|_{z(x)={\rm const}}= 
\f{1}{4\pi i}\int_{s_i^*}  f_j(t) b(x,t) v(t)\;,
\label{varfa}
\ee
where $\{s_i^*\}_{i=1}^{6g-6}$ is the basis of cycles dual to $\{s_i\}_{i=1}^{6g-6}$.
\end{proposition}
\noindent {\bf Proof.}
The proof of (\ref{varfa}) is based on  \cite{KokKor} (Theorem 3).
Namely, the pair $(\Ch,v)$ is a point in the moduli space $\Hc_{\gh}(2,\dots,2)$ 
of Riemann surfaces of genus $\gh$ with an Abelian differential $v$ whose zeros all have multiplicity $2$. Therefore, the space ${\Q_{g}^{0}} $ can be viewed as a subset of $\Hc_{\gh}(2,\dots,2)$ where  $\Ch$ possesses
a holomorphic involution $\mu$, and $v$ is anti-symmetric under this involution.
The proof of (\ref{varfa}) then follows the proof of Lemma 5 of \cite{KorZog2}, to which we refer.

Consider, for example, the derivative with respect to $A_1=\int_{a_1^-}v$, where (see (\ref{abm1})) $a_1^-=\f{1}{2}(a_1-a_1^\mu)$.
According to Th. 3 of \cite{KokKor}, variational formulas for $\vh_i(x)=\fh_i(x)v(x)$ \eqref{basisCH}  on $\Ch$
are given by:
\be
\f{\p\fh_i(x)}{\p (\oint_{a^-_1} v)}=-\f{1}{2\pi i}\oint_{y\in b_1^-} 
\fh_i(y)\bh(x,y) v(y)
\ee
where $\Bh(x,y)=\bh(x,y)v(x)v(y)$ is the  canonical bidifferential on $\Ch$ normalized with respect to the basis (\ref{mainbasis}) on $\Ch$ (requiring that the periods of $\Bh$ along cycles  $a_{j}^+, a^-_{j} $ \eqref{abp} and (\ref{abm1}), \eqref{abm} vanish).
Considering the period of $\oint_{a_k^\mu}v $ as an independent variable, Thm. 3  of \cite{KokKor} implies
\be
\f{\p\fh_i(x)}{\p \left(\oint_{a_k^\mu}v\right)}=-\f{1}{2\pi i}\int_{y\in b_k^\mu} \fh_i(y)\bh(x,y) v(y)\;.
\ee
Recall that  $A_1=\f{1}{2}(\int_{a_1}v-\int_{a_1^\mu}v)$ and $\int_{a_1}v=-\int_{a_1^\mu}v=A_1$. Using the chain rule and taking into account the symmetry
$v(y^\mu)=-v(y)$, we find
\be
\f{\p\fh_i(x)}{\p A_1}=-\f{1}{2\pi i}\int_{y\in b_1} \left(\fh_i(y)\bh(x,y)+ \fh_i(y^\mu)\bh(x,y^\mu)\right)v(y)\;.
\label{fh1}\ee
Along the same lines we also find
\be
\f{\p\fh_i(x^\mu)}{\p A_1}=-\f{1}{2\pi i}\int_{y\in b_1} \left(\fh_i(y)\bh(x^\mu,y)+ \fh_i(y^\mu)\bh(x^\mu,y^\mu)\right)v(y)\;.
\label{fh2}
\ee
Under our choice of canonical basis of cycles in  $H_1(\Ch,\R)$, the set of $a$-cycles is invariant under $\mu$ and therefore
$$
B(x,y)=\Bh(x,y)+\Bh(x,y^\mu) \;,
$$
or, equivalently, 
$$
b(x,y)=\bh(x,y)+\bh(x,y^\mu).
$$
Similarly, 
$f_i(x)=\fh_i(x)+\fh_i(x^\mu)$  (i.e. $v_i(x)=\vh_i(x)+\vh_i(x^\mu)$).

Adding (\ref{fh1}) and (\ref{fh2}) one obtains
$$
\f{\p f_i(x)}{\p A_1}=-\f{1}{2\pi i}\int_{y\in b_1} f_i(y)b(x,y) v(y)\;.
$$
The integrand is anti-symmetric and $b_{1}^-=\f{1}{2}(b_1-b_1^\mu)$, therefore,
$$
\f{\p f_i(x)}{\p A_1}=-\f{1}{2\pi i}\int_{y\in b_1^-} f_i(y)b(x,y) v(y)\;
$$
which implies (\ref{varfa}) when $s_i=a_1^-$, $s_i^*=-2b_1^-$  (recall the normalization $a_i^-\circ b_j^-=\delta_{ij}/2$). Analogously one can verify (\ref{varfa}) for any cycle in the given basis. \QED

Variational formulas for the matrix of $b$-periods and the functions $b(x,y)$ \eqref{bpm} and $\bc(x)$ \eqref{Qx} are derived from the corresponding variational formulas on 
$H_{\gh}(2,\dots,2)$ following the steps in the proof above verbatim.
The resulting formulas are summarized in the following proposition:

\begin{proposition}
Given an arbitrary basis $\{s_i\}_{i=1,\dots 6g-6}$ of $H_-(\Ch, \R)$ and its dual basis $\{s_i^*\}_{i=1,\dots 6g-6}$, the following variational formulas hold
\be
\frac{\p \Omega_{jk}}{\p \P_{s_i}} = \f{1}{2} \int_{s_i^*} f_j f_k v
\label{varO}
\ee
\be
\f{\p b (x,y)}{\p \P_{s_i}} = \f{1}{4\pi i}\int_{ s_i^*} b(x,t)\, b(t,y)v(t),
\label{varB}
\ee
 where $\P_{s_i} := \oint_{s_i} v$ and all derivatives are computed keeping $z(x)$ constant.
Moreover, the variational formula (\ref{varB}) implies the  variational formula for the function $\bc$ (also at $z(x)$ constant):
\be
\frac{\p \bc(x)}{\p \P_{s_i}}=\f{3}{4\pi i} \int_{ s_i^*} h(x,t)v(t)
\label{varQ}.
\ee
\end{proposition}
\section{Symplectic structures on the space of quadratic differentials}
\la{sympquad}
\subsection{Homological symplectic structure on $\Q_g^0$}
The periods of $v$ on $\Ch$ provide a local coordinate system for $\Q_g^0$, and thus we define a Poisson structure on $\Q_g^0$ by the bracket of any pair of periods. 
The ``homological symplectic bracket''  is given by the formula 
\be
\bigg\{
\oint_{s_1} v, \oint_{s_2} v
\bigg\} = s_1\circ s_2, \qquad s_1,s_2 \in H_-(\Ch,\Z).
\la{PSfund}
\ee
By the nondegeneracy of the intersection pairing, the above bracket is also non-degenerate, and hence it defines a symplectic structure.

A set of Darboux coordinates for the homological symplectic structure is obtained by choosing a symplectic basis in $H_-(\C,\R)$; we will therefore use the basis  $\{a_k^-,b_k^-\}_{k=1}^{3g-3}$  (\ref{abm}) to introduce  the homological coordinates
\be
A_k=\int_{a_k^-} v\;,\hskip0.7cm B_k=\int_{b_k^-}v  \ ,\qquad k=1,\dots, 3g-3
\label{homcoord}
\ee
which coincide with  Darboux coordinates up to a  multiplicative factor of   $\sqrt{2}$.
In terms of these coordinates, the  homological symplectic form reads:
\be
{\omega}= 2\sum_{l=1}^{3g-3} dA_l\wedge dB_l
\label{fundam}
\ee
which is manifestly independent of the choice of any basis in $H_-$ as long as the basis of cycles  satisfies \eqref{intH-}.
\subsection{Canonical symplectic structure on $\Q_g$.}
The space $\Q_g$ can be identified with the
cotangent bundle $T^*\Mc_g$ over the moduli space $\Mc_g$ of Riemann surfaces of genus $g$.\footnote{Since the moduli space is an orbifold, the
cotangent bundle is well-defined over the space of Riemann surfaces without automorphisms; the orbifold points are resolved if one considers the
universal covering space of $\Mc_g$ i.e. the Teichm\"uller space $\Tc_g$ and the cotangent bundle over it.}
The complex cotangent bundle $T^*M$ of any complex $n$-dimensional manifold $M$ has a canonical (holomorphic) symplectic structure defined as follows. 
Introduce a system of local coordinates $\{q_i\}_{i=1}^n$ on $M$, and consider the basis $\{\d q_i\}_{i=1}^n$ of the cotangent space. Any cotangent vector 
can be represented as $\sum_{i=1}^n p_i \d q_i$ with some coefficients $p_i$, and the canonical symplectic form on $T^*M$ is defined by
\be
\ocan=\sum_{j=1}^n \d p_j\wedge \d  q_j\;.
\ee
Introducing the {\it symplectic potential} 
\be
\tcan = \sum_{j=1}^n p_j \d q_j\;,
\label{cansp}
\ee
 the symplectic form  can be represented as 
$
\ocan= \d\tcan\;.
$
Neither  $\ocan$ nor $\theta _{{can}}$  depend on the choice of the local coordinates $q_i$, and therefore, they are  globally defined on $T^*M$
\cite{Weinstein}. 

Let us apply this construction when $M$ is the subspace $\Mc^0_g$ consisting of non-orbifold points of $\Mc_g$ (or alternatively assume that $M$ is the Teichm\"uller space $\Tc_g$ with the hyperelliptic locus removed for $g>2$).
According to the Torelli theorem, the period matrix $\O$  uniquely determines the Riemann surface $\CC$.
In fact, in a neighborhood of any non-orbifold  point 
$\CC\in \Mc^0_g$ (in particular not on the hyperelliptic locus for $g>2$), one can choose a subset 
  of $3g-3$ matrix entries of $\O$ to define local analytic
coordinates on $\Mc_g^0$ (see p.216 of \cite{ArbCor}); this subset of entries is denoted by $D$. Thus  
\be
\label{nonhyp}
q_{jk}={\O}_{jk}\;,\hskip0.7cm (jk)\in D\;.
\ee

Any tangent vector in  $T_{[\CC]} \Mc^0_g$ can be identified with a harmonic Beltrami differential $\mu$, and 
the  Ahlfors-Rauch variational formula 
expresses the  variation of $\O_{jk}$ in the direction of $\mu$ as follows:
\be
\label{AhRa}
\langle \mu, d \O_{jk}\rangle=\delta_\mu \B_{jk} =\iint_\CC \mu v_j v_k\;,
\ee
i.e. by the natural pairing of the Beltrami differential $\mu$ with the quadratic differential $v_j v_k$.

The formula \eqref{AhRa} implies that the cotangent vector $\d \O_{jk}\in T^*_{[\CC]}\Mc_g$ can be identified with the quadratic differential $v_j v_k$. In other words, quadratic differentials of this form provide a basis for the space of quadratic differentials as long as $\CC$ is of genus $g\geq 2$ and not hyperelliptic for $g\geq 3$. 
\subsection{Equivalence of the canonical and homological  symplectic structures}
\la{cansym}
It is natural to ask about the relationship between the 2-forms $\o$ and $\o _{{can}}$ on the intersection of the spaces $\Q_g^0$ and $T^*\Mc_g^0$, i.e.
on the space of quadratic differentials with simple zeros on Riemann surfaces without automorphisms.

\begin{theorem}\label{taut}
In the intersection ${\Q_{g}^{0}} \cap T^*\Mc^0_g$ the following identity holds:
\be
\sum_{i=1}^{3g-3} (A_i \d B_i- B_i \d A_i) =\tcan
\label{proptav}
\ee
where $\tcan$ is the canonical symplectic potential (\ref{cansp}) on $T^*\Mc^0_g$. 
Consequently the corresponding  symplectic  2-forms also coincide:
\be
\o=\o _{{can}}\;.
\la{oocan}
\ee
\end{theorem}
\begin{remark}
The canonical 2-form $\o _{{can}}$ can be analytically extended to  $\Tc_g$, where the orbifold points of $\Mc_g$ become
ordinary smooth points; the relation (\ref{oocan}) therefore implies that the form $\o$ can also be extended to  $\Tc_g$ from its original domain of definition i.e. from $\Q_g^0$.
\end{remark}
\noindent 
{\bf Proof.}
Since the hyperelliptic locus is removed from consideration,   an arbitrary quadratic differential from $\Q_g^0$ can be written as a linear combination of $v_j v_k$ , $(jk)\in D$, where $D$ is a set of $3g-3$ matrix entries of the matrix $\B$ yielding independent local coordinates $q_{jk} = \O_{jk}$ on $\M_g$:
\be
Q=\sum_{(jk)\in D} p_{jk} \, v_j v_k\;.
\label{defp}
\ee
 The canonical symplectic form on $T^*\Mc_g^0$
is locally  represented as
$$
\omega _{{can}}=\sum_{(jk)\in D} \d p_{jk}\wedge \d q_{jk} = \d \theta _{{can}}\;.
$$
Since $\qd=v^2$, the definition (\ref{defp}) of $p_{jk}$ can be used to represent the Prym differential $v$ as a linear combination of 
 Prym differentials $v_j v_k/v$ for  $(jk)\in D$:
\be
v=\sum_{(jk)\in D} p_{jk} \f{v_j v_k}{v}\;.
\label{Prymb}
\ee
Integrating this relation over the cycles $a_i^-$ and $b_i^-$ we obtain
\be
A_i=\sum_{(jk)\in D} p_{jk} \int_{a_i^-}\f{v_j v_k}{v}\;\;,\hskip0.7cm
B_i=\sum_{(jk)\in D} p_{jk} \int_{b_i^-}\f{v_j v_k}{v}\;.
\label{ABv}
\ee
Using the variational formulas \eqref{varO} for $\O$, and noticing that for $s_i=a_i^-, b_i^-$ one has $s_i^*=-2b_i^-, 2a_j^-$ respectively, the equations \eqref{ABv} can be further rewritten as
$$
A_i=\sum_{(jk)\in D} p_{jk} \f{\p q_{jk}}{\p B_i}\;,\hskip0.7cm
B_i=-\sum_{(jk)\in D} p_{jk} \f{\p q_{jk}}{\p A_i}\;.
$$
Therefore 
$$
\sum_{i=1}^{3g-3} (A_i \d B_i -B_i \d A_i) = \sum_{(jk)\in D} \sum_{i=1}^{3g-3} p_{jk} \left(\f{\p q_{jk}}{\p B_i} \d B_i+\f{\p q_{jk}}{\p A_i} \d A_i\right)=\sum_{(jk)\in D} p_{jk}\;\d q_{jk} = \theta _{{can}} 
$$
which leads to (\ref{proptav}).
The relation (\ref{oocan}) is obtained by applying the exterior differentiation  to (\ref{proptav}).
\QED

\begin{corollary}
\label{corgen1}
Let $\{q_i\}_{i=1}^{3g-3}$ be any set of local holomorphic coordinates on $\Mc_g^0$, $\{p_{i}\}_{i=1}^{3g-3}$ be the corresponding momenta on $T^*\Mc_g^0$ and $(A_i, B_i)$ be a set of local homological coordinates. The generating function for the symplectic transformation between coordinate systems $(p_{i},q_{i})$ and $(A_i, B_i)$ equals
\be
G^{hom}_{can}=\sum_{i=1}^{3g-3} A_i B_i
\label{gener1}
\ee
\end{corollary}
\noindent {\bf Proof.}
Let 
\be
\theta =2 \sum_{i=1}^{3g-3} A_i \d B_i
\la{sympot}\ee
 be the symplectic potential of $\omega$ in  the coordinates $(A_i, B_i)$; then 
 
\be
\d G^{hom}_{can} = \d \le(\sum_{i=1}^{3g-3} A_i B_i\ri) = \theta-\sum_{i=1}^{3g-3} (A_i \d B_i- B_i \d A_i)
\ee
which  equals      $\theta-\theta_{can}$ according to      Thm. \ref{taut}.
\QED

Although the symplectic potential $\theta _{{can}}$ of the canonical symplectic structure on $T^*\Mc_g$ is independent of the choice of the local coordinates $q_i$ on $\Mc_g$, the symplectic potential $\theta$ (\ref{sympot}) depends on the choice of the homological coordinates. Namely, if $\{A_i,B_i\}$ and  $\{\tilde{A}_i,\tilde{B}_i\}$ are two sets of homological coordinates (related by a constant $Sp(6g-6,\C)$ matrix), then, due to (\ref{proptav}),
\be
\tilde{\theta}-\theta= \d \left(\sum_{i=1}^{3g-3}(A_i B_i-\tilde{A}_i\tilde{B}_i)\right)\;;
\ee
Thus the generating function corresponding to a change of homological Darboux coordinates is 
\be
G^{hom}_{hom}=\sum_{i=1}^{3g-3}(A_i B_i-    \tilde{A}_i\tilde{B}_i)
\la{homhom}
\ee
\begin{remark}
\sf
Theorem \ref{taut} implies that  the entries of $\O$ commute:
\be
\{\Omega_{ij},\Omega_{kl}\}=0\:.
\label{OmOm}
\ee
This relation can also be verified directly using the definition of the Poisson bracket for a pair of functions $\Omega_{ij}$ and $\Omega_{kl}$ along with the 
variational formulas (\ref{varO}). The commutativity of entries in $\O$ can be expressed in terms of 
Riemann bilinear relations for the pair of holomorphic differentials $f_i f_j v$ and $f_k f_l v$. 

\end{remark}

\subsection{The Prym matrix and the change from canonical to homological coordinates}
\label{SecPrym}
The systems of canonical Darboux coordinates $(p_{jk},q_{jk})$ and homological Darboux coordinates $(\sqrt{2}A_i,\sqrt{2}B_i)$ provide distinct polarizations of the space $\Q_g^0$ (i.e. different foliations by 
Lagrangian submanifolds) which we discuss here in more detail.
The Prym matrix (\ref{Prymmat}) plays a key role in this relationship. 

It follows from the definition of the matrix $\Pi$ that the periods of $v$ satisfy the following relation:
\be
B_i=\sum_{j=1}^{3g-3} \Pi_{ij} A_j
\la{BAPi}\ee
This relation holds for the periods of any differential in $H^-$, and in particular for the differentials $v_j v_k/v$:
\be
\int_{b_i^-} \f{v_k v_l}{v} =\sum_{j=1}^{3g-3} \Pi_{ij} \int_{a_j^-} \f{v_k v_l}{v} \;.
\ee
According to the variational formulas (\ref{varO}), the relation can be written as follows:
\be 
 \VF_i [\O]=0\;
\la{vecf}
\ee
where the vector fields $\VF_i$ are defined by
\be
\VF_i=\p_{A_i}+\sum_{j=1}^{3g-3} \Pi_{ij} \p_{B_j}\;,\hskip0.7cm i=1,\dots,3g-3\;.
\la{defvfi}
\ee
Hence the vectors $\VF_j$ span (locally) the tangent bundle of the  vertical foliation in  $T( T^*\M_g)$. 

\begin{proposition}
The vector fields $\VF_i$, $i=1,\dots,3g-3$, given by (\ref{defvfi}), act trivially on functions which depend on the coordinates $\{q_i\}$ only, i.e. they are vertical vector fields for the canonical polarization of $T^*\Mc_g$. Moreover they commute,
 $$[\VF_i, \VF_j]=0\;.$$
\end{proposition}

{\bf Proof.}
The vector fields $\VF_i$ are clearly independent and  act trivially on the canonical coordinates $\{q_i\}$ due to (\ref{vecf}). To prove that $\VF_i$ commutes with $\VF_j$, notice that the
commutator $[\VF_i, \VF_j]$ annihilates any function of the coordinates $\{q_i\}$ which implies (since they span the vertical foliation) that it is a linear combination
of the vector fields $\VF_1,\dots,\VF_{3g-3}$. However, as one sees from the definition (\ref{defvfi}), the vector $[\VF_i, \VF_j]$ is a combination of vectors 
$\p_{B_i}$ only; therefore $[\VF_i, \VF_j]$ must vanish. \QED
\begin{proposition}
The Prym matrix is given by the second Lie derivatives of the function $\f{1}{2}G_{can}^{hom}$ (\ref{gener1}) along vector fields $\VF_i$:
\be
\Pi_{ij}=\f{1}{2}\VF_i\VF_j \left(\sum_{s=1}^{3g-3} A_s B_s\right)\;.
\la{GPrym}
\ee
In addition, the matrix $\Pi$ satisfies the following system of differential equations with respect to the homological coordinates;
\be
\frac {\pa \Pi_{k\ell}}{\pa A_j} +\sum_{r=1}^{3g-3}
 \Pi_{jr} \frac {\pa \Pi_{k\ell}}{\pa B_r} =
 \frac {\pa \Pi_{j\ell}}{\pa A_k} 
+\sum_{r=1}^{3g-3}
 \Pi_{kr} \frac {\pa \Pi_{j\ell}}{\pa B_r}  \ ,\qquad \forall \ j,k,\ell=1,\dots,  3g-3\;.
\la{sysPi}
\ee
\end{proposition}
{\bf Proof.} The expression (\ref{GPrym}) follows directly from the definition of $\{\VF_j\}$; notice that the symmetry of $\Pi$ and the commutativity of the $\VF_i$'s is consistent with \eqref{GPrym}.

Equations (\ref{sysPi}) result from the vanishing of the coefficients in front of $\p_{B_i}$ in the commutator $[\VF_i, \VF_j]$; they are equivalent to the symmetry
of the third Lie derivative $\VF_j\VF_k\VF_l \left(\sum_{s=1}^{3g-3} A_s B_s\right)$ with respect to the interchange of indices $j\leftrightarrow k$.
An alternative proof of 
(\ref{sysPi}) can be obtained using 
variational formulas for the matrix $\Pi$
(see \cite{KorZog2}):
$$
\f{\p \Pi_{kl}}{\p \P_{s_i}} = \f{1}{2} \int_{s_i^*} \f{v_k^- v_l^-}{v}\;.
$$
Using relations between periods of meromorphic second kind differentials ${v_k^- v_l^-}/{v}$ (which have poles of second order at branch points of $\Ch$) and the holomorphic Prym differentials $v_k^-$, the left-hand side of (\ref{sysPi}) can be computed by Riemann bilinear relations 
for the pair of differentials $(v_j^-, \; v_k^- v_l^-/v)$, and it is therefore equal to $2\pi i \sum_{i=1}^{4g-4} \res_{x_i}\left(\f{v_k^- v_l^-}{v}\int^x v_j^-\right)$. This expression is proportional, up to an explicit constant, to
$\sum_{i=1}^{4g-4} (f_k^-)_i (f_l^-)_i (f_j^-)_i$, where $(f_j^-)_i$ are the leading coefficients of the expansion of $v_j^-$ near the branch point $x_i$ in the 
distinguished local parameter $\xi_i(x)=(\int_{x_i}^x v)^{1/3}$: $v_j^-= (f_j^-)_i   d \xi_i+\dots$. Again the $j\leftrightarrow k$ symmetry of this expression for a given $l$ implies (\ref{sysPi}). \QED

The system of equations (\ref{sysPi}) for the Prym matrix seems to be new; together with (\ref{BAPi}) this system contains structural information about the Prym matrix $\Pi$, and therefore we expect it to be useful in the study of the space of quadratic differentials, and in the problem of characterizing
of the corresponding Prymians (a natural analog of the Schottky problem in the case of the moduli space of Riemann surfaces).
\subsection{Computations of homological Poisson brackets}
\label{various_brackets}
For later convenience we now compute some additional Poisson brackets.
\begin{proposition}
\label{propVarO}
For $s\in H_-(\Ch,\R) $ denote by $\P_s=\oint_s v$ the corresponding  period of $v$. The Poisson bracket (\ref{fundam}) implies the following:
\bea
\le\{\O_{jk}, \P_s \ri\}&\&=-\f{1}{2}\oint_s \frac{v_j v_k}{v} \;.
\label{poO}
\\
\le\{f_j(x),\P_s \ri\}&\& = -\f{1}{4\pi i}\oint_{ s}  \frac{v_j(t) B(x,t)}{v(t)v(x)}
\label{pof}
\\
\{ \bc(x), \P_s\}&\&=-\f{3}{4\pi i}\oint_s h(x,t) v(t)
\label{varQ1}
\eea
where  $f_j(x)  =\frac {v_j}{v}$. 
In the Poisson brackets above, the variations are taken assuming that $z(x)$ and $z(y)$ remain fixed under differentiation.
\end{proposition}
\noindent {\bf Proof.}
These expressions follow directly from the definition \eqref{fundam} using the basis of $H_-$ given in \eqref{abm} and the
variational formulas (\ref{varfa}), (\ref{varO}) and (\ref{varQ1}). \QED
\begin{lemma}
The following Poisson brackets hold
\bea
 \{f_i(z), \Omega_{jk}\} &\& = \f{1}{2}d_z\left(f_i\int_{{x_1}}^{{x}}{f_jf_k v} \right)
\label{Of} 
\\
 \{\bc(z), \Omega_{ij}\}&\& =\f{1}{4}(f_i f_j)_{zz}+ \bc f_i f_j + \f{1}{2}\bc_z\int_{x_1}^{{x}}
{f_jf_k v}
\label{OQ} 
\eea
where $z = z(x)$ \eqref{coordz} is kept constant in the computation of the brackets.
\end{lemma}

\noindent {\bf Proof.} Introduce the differentials $V(t)= f_i(t)b(t,x)v(t)$ and $W(t)= f_j f_k v$.
The variational formula (\ref{varfa}) and Riemann Bilinear relations imply that the Poisson bracket $\{f_i(z), \Omega_{jk}\}$ can be expressed as follows:
\be
\{f_i(z), \Omega_{jk}\}=\f{1}{8\pi i}\int_{\p \Ch_0} \left( V(t)\int_{x_1}^t W\right)\;.
\label{Of1} 
\ee
The integral is independent of the choice of base point because the difference between any two choices is the sum of the residues of $V$, which is notoriously zero. For convenience we choose the initial point of integration in the r.h.s. of (\ref{Of1}) to coincide with the zero $x_1$ of $v$.  Assume the boundary $\p  \Ch_0$ is  invariant with respect to the involution $\mu$. 

The differential $W(t)$ has two second order poles: at $t=x$ and $t=x^\mu$. The right-hand side of (\ref{Of1}) is determined via computing the residues at these two points. The differentials $f_j f_k v$ and $v$ are  both skew-symmetric with respect to $\mu$,
therefore the contributions from these two points add, and it's sufficient to compute the residue at $t=x$ which equals
$$
\f{\d}{\d z(x)} \left(f_i(x)\int_{x_1}^x W\right)\;.
$$
This clearly implies (\ref{Of}) where the local coordinate $z$ is used throughout.

To prove (\ref{OQ}) introduce $V(t)= b^2(x,t) v(t)$ and $W= f_i f_j v$. Following the same steps as before and using (\ref{varO}) and (\ref{varQ}), 
\be
\{\bc(z), \Omega_{ij}\}=\f{3}{8\pi i}\int_{\p \Ch_0} \left( V(t)\int_{x_1}^t W\right)\;.
\label{OQ1}
\ee
The integrand has poles on $\Ch$ at $t=x$ and $t=x^\mu$; again residues of this expression at $x$ and $x^\mu$ coincide. The residue at $t=x$ are computed using the asymptotics (\ref{asymph}) of $b^2(x,t)$ as $t\to x$. The result, multiplied by 2 to take into account contributions of both $x$ and $x^\mu$, is the r.h.s. of (\ref{OQ}).
\QED

\begin{remark}
It is instructive to verify  the Jacobi identity for the Poisson bracket of various triples of functions directly. The simplest, but rather non-trivial verification, is the Jacobi identity for the triple 
$\left(\Omega_{ij}\,,\;  \Omega_{kl}\,, \;q(z)\right)$. All entries of the matrix of $b$-periods Poisson-commute (\ref{OmOm}), and the Jacobi identity is equivalent to the symmetry of the expression $\{\{q(z),\Omega_{ij}\},\Omega_{kl}\}$ under  the interchange of  $\Omega_{ij}$ and $\Omega_{kl}$. This symmetry can be verified by a straightforward computation using (\ref{Of}) and (\ref{OQ}).
\end{remark}

\begin{proposition}
\label{Poissonq}
Assuming that the  coordinates $z$ and $\zeta$ remain independent of the homological coordinates on $\Q_g^0$ one has
\bea
\f{4\pi i}{3}\{\bc(z), \bc(\zeta)\}=&\& \f{1}{2}h_{zz}(z,\zeta)- \f{1}{2} h_{\zeta\zeta}(z,\zeta) +2 h(z,\zeta)(\bc(z)-\bc(\zeta))
\nonumber \\
&\&+\bc_{z}(z)\int_{0}^z h(\zeta,t)dt-\bc_{\zeta}(\zeta)\int_{0}^\zeta h(z,t)dt
\label{QQ}
\eea
where $h(z,\zeta)=b^2(z,\zeta)$.
\end{proposition}
\noindent {\bf Proof.} The proof is parallel to the proof of (\ref{OQ}); the contributions to the integral over the boundary of the fundamental domain come from residues at $x,x^{\mu},y$ and $y^\mu$.
\QED

\noindent

In terms of  $u(z)$ \eqref{uq}, Proposition \ref{Poissonq} takes the form:
\bea
\f{4\pi i}{3}\{u(z), u(\zeta)\}
=&\&\left(\f{1}{2}\h^{(\zeta)}_{zzz}-2 \h^{(\zeta)}_{z}u(z) - \h^{(\zeta)}(z) u_z(z)\right)
-\left(\f{1}{2}\h^{(z)}_{\zeta\zeta\zeta}-2 \h^{(z)}_{\zeta}u(\zeta) - h^{(z)}(\zeta) u_\zeta(\zeta)\right)\;.
\label{Poissonu}
\eea
The function $h^{(z)}(\zeta)$ is given by \eqref{hxy}, and \eqref{Poissonu} is re-expressed as \eqref{QQ} using the property $(\h^{(z)}(\zeta))_\zeta=(\h^{(\zeta)}(z))_z=\h(z,\zeta)$ which follows directly from the definition of $\h^{(z)}(\zeta)$.

\subsection{Commuting homological flows}      
   \label{secthomo}
A polarization of the symplectic structure also defines a maximal set of commuting Liouville integrable Hamiltonian systems 
  on the Torelli cover, $\wt {\Q_g^0}, $ of $\Q_g^0\subset T^*\Mc_g$.
  Following section \ref{cansym} we  choose a subset 
  $D$ of $3g-3$ entries of the period matrix $\Omega$ as canonical Darboux coordinates $q_{jk}$; the corresponding momenta are  the coefficients $p_{jk}$ of a holomorphic quadratic differential in the basis
   $\{v_j v_k\}$ for $(jk)\in D$.

   The following family of Hamiltonians $H_1,\dots,H_{3g-3}$ commute with respect to the canonical Poisson structure and define an integrable system on $\Q_g^0$:
   \be
   H_i=\f{1}{2\pi} A_i^2\;,\hskip0.7cm i=1,\dots,3g-3
   \la{HA}
   \ee
   where the homological coordinates $A_1,\dots,A_{3g-3}$ are the periods of $v$ over the cycles $a_j^-$ in the symplectic basis for $H_-$ chosen  in (\ref{abm}). The corresponding canonical basis of cycles on $\Ch$ is shown in Fig.\ref{canbasCh}.
   Denote  by $t_i$ the time variables for the Hamiltonians $H_i$. The choice of quadratic expressions (\ref{HA}) is justified below as they are exactly the action variables of this  integrable system; moreover, $A_i^2$, in contrast to $A_i$, does not depend on the choice of the sign of $v=\sqrt{\qd}$.
    
 According to (\ref{poO}) the equations of motion for $\O_{jk}$ take the form
   \be
   \f{d \O_{jk}}{d t_i}=\{H_i, \O_{jk} \}=\f{1}{2\pi}A_i\int_{a_i^-}\f{v_j v_k}{v}.
   \la{eqmoq}
   \ee
The evolution of   the $(3g-3)\times (3g-3)$ Prym matrix $\Pi$ (\ref{Prymmat}) is given by

  \be
  \f{d \Pi_{jk}}{d t_i}=\{H_i, \Pi_{jk}\}=\f{1}{2\pi}A_i\int_{a_i^-}\f{v^-_j v^-_k}{v}
  \la{eqmoPrym}
  \ee
  where $v_j^-$ are the Prym differentials (\ref{Prym}).
  According to Torelli's theorem,  the pair of matrices $(\O_{jk},\Pi_{jk})$
determines the moduli of the canonical cover $\Ch$,  and therefore the equations of motion (\ref{eqmoq}) and (\ref{eqmoPrym}) completely define
the dynamics of a pair $(\CC,\qd)$ up to multiplication of $\qd$ by a constant. This constant is determined by requiring that the Hamiltonian (\ref{HA}) remains constant under the $t_i$-evolution.

  The set of commuting Hamiltonians $\{H_1,\dots,H_{3g-3}\}$ can be split into two groups: the Hamiltonians 
  \be\la{hambase}
  H_j=\f{1}{8\pi }\left(\int_{a_j-a_j^\mu} v\right)^2\;,\hskip0.7cm j=1,\dots,g
  \ee
are expressed via   periods of $v$ over of cycles in $H_-$ arising from homologically non-trivial cycles on $\CC$, and
  \be
  H_{g+m}=\f{1}{4\pi}\left(\int_{\tilde{a}_m} v\right)^2\;,\hskip0.7cm m=1,\dots,2g-3
  \la{hamcov}\ee
 are obtained by  integrating $v$ around branch points of $\Ch$ (see Figure \ref{genpi1}).

 Denote the level sets of the Hamiltonians, $H_i$, in $\Q_g^0$ by $\Q_H$. In the general formalism (see \cite{Woodhouse}), the action variables are the periods of the canonical symplectic potential $\theta _{{can}}$ over a set of closed loops $\varpi_1,\dots,\varpi_s$ in $\pi_1(\Q_H)$:
  \be
  I_\ell= \f{1}{2\pi}\int_{\varpi_\ell} \theta _{{can}}
  \la{defact}
  \ee
  Since $\Q_H$ is a Lagrangian submanifold, the restriction of the form $\theta _{{can}}$ on $\Q_H$ is closed, and the actions only depend on the level-set and the homotopy class of the contour in $\Q_H$.
  The Dehn twist on (both copies of) $\CC$ along the corresponding generator $\a_s$ of  $\pi_1(\CC,x_0)$ gives a natural choice of paths in the space $\Q_H$;
   these paths will be denoted by $\varpi_s$ for $s=1,\dots,g$. Each Dehn twist can be performed while keeping all $a^-$-periods of $v$ fixed, and thus the deformation path lies inside of $\Q_H$. Choose the direction of the Dehn twist such that under the deformation along $\varpi_s$ the corresponding period $B_s$ is increased by $A_s$  according to the Picard-Lefshetz formula
  \be
  B_s\to B_s+A_s
  \la{deforB}
  \ee
  while all other $B_i$'s remain constant.
  For the second set (\ref{hamcov}) of commuting Hamiltonians, the deformation path $\varpi_{g+k}$ in $\pi_1(\Q_H)$ is chosen so as to generate the interchange of the zeros $x_{2k+1}$ and $x_{2k+2}$. This deformation is  half of a  Dehn twist on $\Ch$ along the path encircling the branch cut $[x_{2k+1},x_{2k+2}]$ on $\Ch$ (in homology such a path corresponds to the cycle $a^-_{g+k}$).
  The corresponding period $B_{g+k}$ also transforms via (\ref{deforB}) under this deformation.
  
  One can now integrate the canonical symplectic potential $\theta _{{can}}$ along $\varpi_s$ using (\ref{proptav}):
  $$
  \oint_{\varpi_s}\!\!\! \theta _{{can}}= \oint_{\varpi_s} \sum_{i=1}^{3g-3} (A_i d B_i-B_i dA_i)=\sum_{i=1}^{3g-3} A_i \int_{\varpi_s}dB_i\;.
 $$
   According to (\ref{deforB}), $\int_{\varpi_s}dB_i= A_s\delta_{is}$. The action-angle variables are therefore given by
   \be
   I_s=\f{1}{2\pi } A_s^2\ ,\qquad \varphi_s=2\pi \f{B_s}{A_s}\;, 
   \la{actions}
   \ee
   and thus the action variable $I_s$ coincides with the Hamiltonian $H_s$.
   The canonical symplectic form then has the standard expression
   \be
   \o _{{can}}=2\sum_{i=1}^{3g-3} d A_i\wedge dB_i=\sum_{i=1}^{3g-3} d I_i\wedge d\varphi_i\;.
   \ee

   Therefore the standard form of the flows generated by the Hamiltonians $H_i=I_i$ in action-angle variables is
   \be
   I_s=const\;, \hskip0.7cm
   \varphi_i=\varphi_i^0+  t_i\;, \hskip0.7cm 
   \varphi_s=const,\;\;\;\; s\neq i
   \ee
    which can be expressed in homological coordinates $(A_s, B_s)$ as follows
   \be
   A_s=const\;,\hskip0.7cm
   B_i=B_i^0 + \f{1}{2\pi}t_i A_i\;, \hskip0.7cm
   B_s=const,\;\;\;\; s\neq i\;.
   \ee
The flows with respect to the times $t_1,\dots,t_g$ can be interpreted geometrically as  shear flows around (the lift to $\Q_g^0$ of) the  boundary component $\delta_0$ of the Deligne-Mumford compactification.
   The flows with respect to the  times $t_{g+1},\dots,t_{3g-3}$ can be interpreted as shear flows around the component $D_{deg}$ of the boundary of $\Q_g^0$ which
   consists of quadratic differentials with one double and $4g-2$ simple zeros \cite{KorZog2}.
\section{Symplectic structures on the space of projective connections}
\label{SecCan}
Denote by $\Proj_g$ the space of pairs $(\CC,\proj)$ where $\proj$ is a holomorphic projective connection on the Riemann surface $\CC$ of genus $g$.
The space $\Proj_g$ is an affine bundle over the moduli space $\Mc_g$ of dimension $6g-6$.

Given a projective connection $S_0$ on each Riemann surface $\CC\in \Mc_g$ which depends holomorphically on the point of $\Mc_g$, any other projective connection $S$ on $\CC$ can be decomposed into the sum of $S_0$ and a holomorphic quadratic differential $2\qd$ on $\CC$:
\be
S=S_0+2\qd
\la{SS0}
\ee
where the factor of 2 is introduced for convenience. Therefore any choice of the reference projective connection $S_0$ (which is locally defined on $\Mc_g$
but globally may depend  on a marking of $\CC$) defines an isomorphism
$$
F^{(S_0)}\;:\;\Q_g\to \Proj_g\;.
$$
Using the isomorphism  $F^{(S_0)}$ we induce a  symplectic form $\o_0$  on $\Proj_g$  from the canonical symplectic form on $\Q_g$. 
If we choose a different holomorphically varying reference projective connection $S_1$ on $\CC$, the map $F^{(S_1)}$ induces another symplectic form $\o_1$ on $\Proj_g$ from the canonical symplectic structure on $T^*\Mc_g$. In general there is no reason to expect the symplectic forms $\o_0$ and $\o_1$ 
to  coincide (see Remark \ref{noncon} below).
 
 \begin{definition}
 \label{defequiv}
 Two  holomorphically varying  projective connections $S_0$ and $S_1$ are  {\it equivalent} if $\o_0=\o_1$, namely, they induce the same symplectic structure on $\Proj_g$.
 \end{definition}

The  equivalence of two projective connections is characterized as follows. Given a local coordinate system $\{q_i\}$ on $\Mc^0_g$, the differentials $\d q_j$ are identified with a basis of the space of holomorphic quadratic differentials. Therefore the quadratic differentials $S-S_0$ and $S-S_1$ can be expressed as  linear combinations of $\{\d q_j\}$:
\be
\proj -S_0=2\sum_{i=1}^{3g-3} p^0_i \d q_i\  \qquad \proj -{S_1}=2 \sum_{i=1}^{3g-3} p^1_i \d q_i\;.
\ee

Denote the corresponding symplectic potentials by $\theta_0 = \sum_{i=1}^{3g-3} p^0_i \d q_i$, $ \theta_1 = \sum_{i=1}^{3g-3} p^1_i \d q_i$.   The corresponding symplectic forms coincide if and only if $\d(\theta_0-\theta_1)=0$ which is equivalent to the existence of a local holomorphic function $G_0^1$ on the base moduli space such that 
$$
\theta_0-\theta_1 = \d G_0^1 \;.
$$
Geometrically,   $S_0$ and $S_1$ within the same equivalence class determine two Lagrangian embeddings of the base moduli space 
into $\Proj_g$. The function $G_0^1$ is the corresponding generating function in the sense of symplectic geometry.

\subsection{Equivalence of  Bergman projective connections for different markings} 
\la{BCdm}

The definition of the reference projective connection $S_0$ may involve a marking  of $\CC$.
For example to define the Bergman projective connection, we need  a Torelli marking of $\CC$. If  $S_0$ is chosen to be  the Schottky projective connection, $\CC$ must be Schottky marked.
 Finally, if $S_0$ is
a Bers projective connection (used by Kawai in \cite{Kawai}), then $\CC$ must be  Teichm\"uller marked. 
Therefore given two (locally) equivalent projective connections, the corresponding generating function  depends on the two markings.

\begin{proposition}
\label{propGstt}  
Let  $\sigma$ be an $Sp(2g,\Z)$ matrix (\ref{symint}).
The Bergman projective connections $S_B$ and $S_B^\s$ (given by (\ref{transSB})) are equivalent in the sense of Definition \ref{defequiv}.
The generating function of the change between the corresponding Lagrangian embeddings of the base is
 \be
G_{B}^{\sigma}  = 6\pi i \log  {\rm det} (C\Omega+D)\;.
\la{GBs}
\ee
\end{proposition}
 {\bf Proof.}  
Due to (\ref{transSB}) the difference of symplectic potentials corresponding to $S_B$ and $S_B^\sigma$ on $\Projm_g$ is given by:
\be
\theta_B - \theta^{\s}_B = 
6\pi i \sum_{ 1\leq j\leq k\leq g} v_jv_k \frac \pa {\pa \O_{jk}} \ln \det (C \O+D)\;.
\label{QQs}
\ee
The Ahlfors-Rauch formula (\ref{AhRa}) together with the chain rule allows one to express the formula (\ref{QQs}) as follows:
\be
\theta_B-\theta_B^{\s}  = 
6\pi i \d \log {\rm det} (C\Omega+D)
\label{phiphi}
\ee 
which is equivalent to (\ref{GBs}).
\QED

\begin{remark}\la{noncon}
The statement that symplectic forms corresponding to different Torelli markings are equal is rather non-trivial when expressed in terms of homological coordinates. The homological coordinates $(A_i,B_i)$, corresponding to the quadratic differential $\qd$, and the homological coordinates $(A^\s_i,B^\s_i)$, corresponding to the quadratic differential $\qd^\s=\qd+1/2(S_B-S_B^\s)$, are related in a highly non-trivial way because the corresponding canonical covers $\Ch$ and $\Ch^\s$ are different. Nevertheless, Proposition \ref{propGstt} implies the relation $\sum_{i=1}^{3g-3} \d A_i\wedge \d B_i=\sum_{i=1}^{3g-3} \d A^\s_i\wedge \d B^\s_i$
which would be  hard to check directly.
\end{remark}

\subsection{Wirtinger and Schottky projective connections.}
\label{Lagrchange}

In this section we show that  the  Schottky and Wirtinger projective connections introduced in Section \ref{Cancov} belong to the same 
equivalence class (Def. \ref{defequiv}) as the Bergman projective connection and find the corresponding generating functions. 
The current proof of equivalence between Bergman and Bers 
projective connections  requires a comparison of our results with the results of Kawai \cite{Kawai} showing that the monodromy mapping is a local symplectomorphism; it is postponed until Sect. \ref{SectBers}.

 \paragraph{Wirtinger projective connection.}
\begin{proposition}
\label{propWirt}
The Wirtinger  and  Bergman projective connections are equivalent.  The generating function of the change between the corresponding Lagrangian embeddings is
\be
G_{B}^{W} = -\f{24\time 4\pi i }{2^g+4^g} \log \left(\prod_{\beta\;\; even}\theta[\beta](0)\right)\;.
\la{GBW}
\ee
\end{proposition} 
{\bf Proof.} As a corollary to the definition of the Wirtinger projective connection (\ref{Swirt}),
the  difference between 
 $\theta_B$ and $\theta_{W}$ 
can be expressed as a total derivative,  
\be
 \theta_B - \theta_{W}  = -\f{24\time 4\pi i }{2^g+4^g}\d \log \left(\prod_{\beta\;\; even}\theta[\beta](0)\right)\;,
\la{SBSW}
\ee
which immediately implies (\ref{GBW}) following the same steps as in Prop. \ref{propGstt} .
\QED

The generating function (\ref{GBW}) is singular on the divisor where some theta-constants vanish, i.e. the divisor where the 
Wirtinger projective connection is singular.

\paragraph{Schottky Projective connection.}
Recall that any nontrivial element $\g$ of the Schottky group $\Gamma_{Sch}$ is loxodromic, and it is characterized by fixed-points $a_\g,b_\g$ and multiplier 
$q_\g$ with $0<|q_\g|<1$. The transformation $w\to \g w$ is then defined by the equation
\be
\f{\gamma w -a_\g}{\g w - b_{\g}}=q_\g \f{ w -a_\g}{ w - b_\g}\;.
\ee

 The Bowen-Zograf ``F-function" is defined  on the Schottky space \cite{Bowen,Zograf1990} (see also \cite{MacTak}) by the absolutely convergent series
\be
F=\prod_{{\gamma}}\prod_{m=0}^\infty (1-q_\gamma^{1+m})
\label{defF}
\ee
where $\g$ runs over all distinct primitive conjugacy classes in $\Gamma_{Sch}$ excluding the identity.

The function $F$ is in fact holomorphic on the Schottky space and was first briefly introduced by R. Bowen  in 1979 \cite{Bowen}. Later it was rediscovered and extensively studied by P. Zograf \cite{Zograf1990} in the context of the holomorphic factorization of the determinant of the Laplace operator on $\CC$.
The characteristic property of $F$ is that under an infinitesimal deformation of the conformal structure by a Beltrami differential $\mu$, it
satisfies the following equation:
\be
\delta_{\mu}\log F=-\f{1}{12\pi i}\iint_{C}(S_{B}-S_{Sch})\mu
\label{eqF}
\ee
where $S_B$ is the Bergman projective connection corresponding to the canonical bidifferential normalized along the system of generators 
defining the Schottky group. This fact can be deduced from \cite{Zograf1990} and \cite{MacTak}, and while some details of the computation are missing in these papers, they can be filled in with some effort.

\begin{proposition}
\label{propSchot} The Schottky and Bergman projective connections are equivalent.
 The generating function of the change  between Lagrangian embeddings is  
 \be
G_{B}^{Sch}=-6\pi i \log\,F\;.
\la{GBS}
\ee
\end{proposition}
{\bf Proof.}
Denote the symplectic potential on $\Proj_g$ corresponding to the Schottky projective connection by $\theta_{Sch}$.
Due to (\ref{eqF}) the difference of symplectic potentials $\theta_{B}$ and $\theta_{Sch}$ is
 \be
\theta_{B}-\theta_{Sch} =  -6\pi i\d \log\,F
\label{ththc}
\ee
where $F$ is the Bowen-Zograf function $F$ (\ref{defF}). This  implies (\ref{GBS}).
\QED

\begin{remark} 
 These examples may give  the impression that any reference projective connection belongs to the same equivalence class as long as it depends  holomorphically on the  moduli. This is false as the following simple counterexample shows. Choose $S_0 = S_{B} + \Omega_{11} v_2^2$. The exterior derivative of
 $S_0-S_B=\Omega_{11} v_2^2 = \Omega_{11} \d \O_{22}$ equals $\d \O_{11} \wedge \d \O_{22}$ which is nowhere vanishing on the moduli space since $v_1^2$ and $v_2^2$ are linearly independent. Therefore although both  $S_B$ and $S_0$  depend holomorphically on the moduli, they correspond to different symplectic structures on $\Proj_g$.
\end{remark}
\subsection{Covering of $\Proj_g$ by charts of homological coordinates}
In this section we show that for any point $S\in \Proj_g$ there exists a Torelli marking $\tau$ such that the quadratic differential $S-S_B^{\tau}$ has only simple zeros. Thus the space $\Proj_g$, can be covered by charts from  homological coordinates.
The statement seems intuitively obvious, but  its rigorous proof is not completely trivial.
 
\begin{lemma}
\label{charts} Let $Q$ be a quadratic differential with higher order zeros, i.e. $Q\in Q_g\setminus\Q_g^0$. There is a symplectic transformation \eqref{symint} $\s$ such that the differential $Q^\s=Q+\f{1}{2}(S_B-S_B^\sigma)$ has only simple zeros .
\end{lemma}
\noindent{\bf Proof.}
It will be sufficient to consider symplectic matrices of the form $\s = \le(\begin{array}{cc}
D & \1 \\
D-\1 & \1
\end{array}\ri) \in Sp_g(\Z)$ where $D\in GL_g(\Z)$ is {\it diagonal}. 
\be
\qd^\s =  \qd +6\pi i {\bf v}(x) (\Omega+ D)^{-1}  {\bf v}^t(x)
\ee
Let $\Delta:= \inf_{j} |D_{jj}|$, and if the matrix $D$ has $\Delta$ sufficiently large, we have
\be
\qd^\s =  \qd +6\pi i \overbrace{ \sum_{j=1}^{g} \frac {v_j^2}{D_{jj}}}^{\mathcal O(\Delta^{-1})} + 6\pi i \overbrace{{\bf v}(x) D^{-1}\O (\Omega+ D)^{-1}  {\bf v}^t(x)}^{=\mathcal O(\Delta^{-2})}
=\qd + V_D +W_D 
\ee
The estimate $\mathcal O(\Delta^{-1})$ is the $\sup$ norm obtained by taking the supremum of the absolute values of the trivialization in a fixed, finite cover of $\CC$ by local coordinates.
By taking $\Delta$ sufficiently large, the term $W_D$ can be neglected, and using the well known fact that the holomorphic differentials $v_j$ cannot have a common zero for all $j=1,\dots,g$, one can show that the term $V_D$ can also be chosen to have only simple zeros .
Therefore, $V_{nD}+W_{nD}$ for $n\in \mathbb N$ large enough, is a sufficiently generic perturbation of $\qd$ so that any zero of multiplicity $k\geq 2$ splits into $k$ simple zeros  in a neighborhood of a multiple zero of $\qd$. \QED
\section{Schwarzian and linear second order equations on Riemann surfaces}
\label{Seceq}
Let $\proj$ be a holomorphic projective connection on a Riemann surface $\CC$, and consider a linear equation of second order
\be
\varphi'' + \f{\proj}{2}\varphi=0\;.
\label{Sch1}
\ee
where prime denotes the holomorphic derivative with respect to any local coordinate. The equation \eqref{Sch1} is invariant under a change of the local coordinate  if
$\varphi$ transforms locally as a $-1/2$-differential \cite{HawSch}.

The ratio $f=\varphi_1/\varphi_2$ of two linearly independent solutions of (\ref{Sch1}) solves the Schwarzian equation 
\be
{\mathcal S}(f,\xi) = S(\xi)
\la{Schwarz}\ee
where $\xi$ is a local coordinate on $\CC$.

The projective connection $S$ can be represented as a sum of the Bergman projective connection corresponding to some choice of Torelli marking on $\CC$ and a quadratic differential. The differential equation (\ref{Sch1}) then takes the form
\be
\label{Schrint}
\varphi'' +\frac{1}{2} \le( S_B+2 \qd\ri)\varphi =0 \;.
\ee 

Assume the Torelli marking is such that the quadratic differential $\qd$ has only simple zeros, and  introduce the canonical cover $\Ch$ by
the equation $v^2=\qd$ (see Sect. \ref{Cancov}).  The zeros of the Abelian differential $v$ on $\Ch$ all have multiplicity $2$, and therefore there exists a section (unique up to a sign) $\spin$ of a spin 
line bundle over $\Ch$ such that $\spin^2=v$. 
Let us now define a function of the coordinate $z(x)=\int_{x_1}^x v$,
$$
\psi(z(x)):=\phi(x) \spin(x)\;.
$$

If one chooses the generators of the fundamental group of $\pi_1(\CC\setminus\{x_i\}_{i=1}^{4g-4},x_0)$ according to (\ref{geninS2}), the function $\psi(z)$ is single-valued on the "first"  sheet of the double cover of the fundamental domain of $\CC$ shown in Fig.\ref{coverfund} outside of the branch points $x_i$;
the coordinate $z$ is also well-defined in the same domain.

The equation (\ref{Schrint}) can  be rewritten in a coordinate-independent way in terms of the function $\psi$ as follows:
\be
d\left(\f{d\psi}{v}\right)+\left(\f{S_B-S_v}{2v}+v\right)\psi=0\;,
\label{inv2}
\ee
and in terms of the coordinate $z$ in the form
$$
\psi_{zz} + \le(\bc  +1\ri) \psi=0\;
$$
where the function $\bc$ is given by (\ref{Qx}).
Introducing the  function $u(z)=-(\bc(z)+1)$ \eqref{uq} this equation  takes the particularly simple form
\be
\psi_{zz} -u \psi=0\;.
\label{Sch2}
\ee
\begin{remark}
The equation (\ref{Sch2}) has apparent singularities at the zeros of $v$. Indeed, the differential $v$ has a second order zeros on $\Ch$
at $x_i$, and therefore $v\sim \hat{\sigma}^2 d\hat{\sigma}$ where $\hat{\sigma}$ is a local parameter on $\Ch$ near $x_i$. Then $z-z(x_i)\sim \hat{\sigma}^3$ as $x\to x_i$;
thus $\hat{\s}\sim (z-z(x_i))^{1/3}$ as $x\to x_i$, and for the local parameter $\s=\hat{\s}^2$ near $x_i$ on $\CC$ one has $z-z(x_i)\sim \sigma^{3/2}$.
Therefore, $\d z/\d\s\sim \s^{1/2}\sim (z-z(x_i))^{1/3}$ and
$$
\psi=\phi\sqrt{\d z}\sim (\phi\sqrt{\d\s})(z-z(x_i))^{1/6}\;.
$$
The function $\phi\sqrt{\d\s}$ is regular at the zeros of $Q$ and $\psi''/\psi\sim -5/(36z^2)$. Thus the potential $u$ 
 has second order poles at the zeros of $v$ corresponding to $z=z(x_i)$ with the local behavior
$$
u(z)=-\frac{5}{36 (z-z(x_i))^2}+ \mathcal O(1)\;.
$$
These singularities are artefacts of the choice of coordinate $z$. In the original form of the equation (\ref{Sch1}), the zeros of $\qd$ are ordinary regular points.
\end{remark}
Given two linearly independent solutions $\psi_{1,2}$ of (\ref{Sch2}),
their Wronskian matrix $\Psi$ solves the first-order ODE which can be written in either invariant form or in terms of the coordinate $z$ as 
\be
d\Psi= \left(\ba{cc} 0 & v \\
uv & 0  \ea \right)\Psi\; \ \ \Leftrightarrow \ \ 
\label{matrix1}
 \f{d\Psi}{d z}= \left(\ba{cc} 0 & 1 \\
u(z) & 0  \ea \right)\Psi\;.
\ee
We assume that the matrix $\Psi$ satisfies the initial condition
\be
\label{IVP}
\Psi(z_0) = \1\ , \ \ \ z_0 = z(x_0), 
\ee
where $x_0$ is the corner of the fundamental polygon on the "first" copy of the fundamental domain of $\CC$ shown in Fig.\ref{genpi1}.

The following notations for a basis in $sl(2,\C)$ is used:
$$
\sigma_-= \left(\ba{cc} 0 & 0 \\
1 & 0  \ea \right)\;\;,\hskip0.7cm
\sigma_+= \left(\ba{cc} 0 & 1 \\
0 & 0  \ea \right)\;\;,\hskip0.7cm
\sigma_3= \left(\ba{cc} 1 & 0 \\
0 & -1  \ea \right)\;.
$$
Starting from the solution $\Psi$ of (\ref{matrix1}), introduce the following bilinear expressions:
\be
\Ld(x)=\Ldm (x)=\Psi^{-1}(x) \sigma_- \Psi(x)
\label{lambdam}
\ee
\be
\Ldp(x)=\Psi^{-1}(x) \sigma_+ \Psi(x)\;\;\;\hskip0.7cm
\Ldt(x)=\Psi^{-1}(x) \sigma_3 \Psi(x)
\label{lambdap3}\;.
\ee
The following simple lemma, whose proof is elementary using \eqref{matrix1}, will be an important technical tool.
\begin{lemma}
The functions $\Ld(z)$, $\Ldp$ and $\Ldt$ satisfy the following equations:
\begin{itemize}
\item[--]
The third order equation for $\Ld$:
\be
\Ld_{zzz}-4 u(z) \Ld_z-2 u_z \Ld=0\;.
\label{3rdor}
\ee
\item[--]
The following two equations express $\Ldp$ and $\Ldt$ via $\Ld$ and its derivatives
\be
\Ldt= -\Ld_z \;,
\label{derF}
\ee
\be
\Ldp=-\f{1}{2} \Ld_{zz}+u(z)\Ld\;.
\label{derFtt}
\ee
\end{itemize}
\end{lemma}

The third order equation (\ref{3rdor}) can also be written in the following forms
\be
(\Ld_{zz}  - 4u\Ld)_z = -2u_z \Ld\;,
\label{fullder}
\ee
\be
\Ld_{zzz}=  2u_z\Ld-4u(z)\Ldt\;.
\label{derFttt}
\ee

\subsection{Transition and monodromy matrices: monodromy representation}
\la{secmon}

Define the transition matrix
\be
T(z_1,z_2)=\Psi(z_1)\Psi^{-1}(z_2)\;,
\label{transdef}
\ee
which is  independent of the normalization point $x_0$ of the solution $\Psi$.
The monodromy matrices can be defined in terms of transition matrices as follows:
\be
M_\g= T(z_0+\P_\g,z_0)
\label{MTRAN}
\ee
where $z_0=\int_{x_1}^{x_0}v$.
The notation for the transition matrix $T(z_0+\P_\g,\,z_0)$ explicitly refers to  the values of the
$z$-coordinate at the initial point and at the endpoint; it is  assumed that the transition matrix between $z_0$ and $z_0+\int_\g v$ is computed along the path $\g$.

 Changing the initial point $z_0$ to $\tilde{z}_0$,  
the corresponding monodromy (anti)representation $\tilde{M}_\g$ transforms by  conjugation;
\be
\tilde{M}_\g= T(\tilde{z}_0+\P_\g,\tilde{z}_0)=
 T(\tilde{z}_0+\P_\g,z_0+\P_\g) M_\g  
T(z_0,\tilde{z}_0)=T^{-1}(z_0,\tilde{z}_0) M_\g T(z_0,\tilde{z}_0)\;.
\ee
Therefore, $\tr M_\g=\tr \tilde{M}_\g$, and  $\tr M_\g$ is independent of the 
initial point $z_0$.

\paragraph{$SL(2,\C)$ monodromy representation of the linear system (\ref{matrix1}).}

The monodromy representation of the linear system (\ref{matrix1}) is closely related to the monodromy representation of the corresponding Schwarzian equation (\ref{Schwarz}). A solution $f$ of the  Schwarzian equation  (\ref{Schwarz})
transforms as follows along a closed loop $\g$:
\be
f\to \frac{a f+ c}{b f + d}\ ,
\la{ffg}\ee
and therefore  the matrix of coefficients 
\be
M_\g=\left(\ba{cc} a & b\\ c & d\ea\right) 
\la{MgS}\ee
is defined only up to an overall sign.
 Thus the monodromy matrices of the Schwarzian equation only define a $PSL(2,\C)$
anti-representation of $\pi_1(\CC,x_0)$ (for convenience $M_\g$ is defined in (\ref{MgS}) by the transposition of the matrix of the 
M\"obius transformation \eqref{ffg}).

On the other hand, one can assign well-defined $SL(2,\C)$ monodromy matrices to the matrix equation 
 (\ref{matrix1}). Namely, the 
coefficients $v$ and $uv$
of (\ref{matrix1}) are Abelian differentials on the canonical cover $\Ch$ which are skew-symmetric under the canonical involution.
However, Lemma \ref{S2lem} allows us to define $SL(2,\C)$ monodromy matrices of equation (\ref{matrix1})  as we  choose the generators $\{\a_i,\b_i,\g_i\}$ of  $\pi_1(\CC\setminus\{x_i\}_{i=1}^{4g-4},x_0)$ satisfying \eqref{geninS2}. Under this choice the matrix of coefficients in \eqref{matrix1} does not change along the loops $\a_i$ or $\b_i$.
The lemma below guarantees this construction gives a well defined  $SL(2,\C)$  monodromy representation of $\pi_1(\CC,x_0)$.

\begin{lemma} 
Choosing generators of the fundamental group $\pi_1(\CC\setminus\{x_i\}_{i=1}^{4g-4},x_0)$ which satisfy Lemma \ref{S2lem} allows to define an $SL(2,\C)$ monodromy (anti)representation of $\pi_1(\CC,x_0)$ corresponding to the matrix equation (\ref{matrix1}). This lifts the $PSL(2,\C)$ monodromy (anti)representation of
the corresponding Schwarzian equation to $SL(2,\C)$.
\end{lemma}
{\bf Proof.}
Let $\xi$ be a local coordinate on $\CC$ around $x_j$ such that $\qd = \xi \d \xi^2$,  and rewrite \eqref{Schrint} in this coordinate.  The coefficients are locally holomorphic in $\xi$, and hence any fundamental local matrix solution $\wt \Psi$  is single valued around $\z=0$. The matrix $\Psi$ is then of the form $\Psi(\xi) =\xi^\frac 1 4  \wt \Psi(\xi) C$, for some invertible constant matrix $C$. This shows that $\Psi$ has  local monodromy $e^{\pi i/2}=i$, and since this is a scalar, the statement holds independently of the chosen local solution $\wt \Psi$. We now choose the set of generators $\{\a_j,\b_j,\g_k\}$ 
of $\pi_1(\CC\setminus\{x_i\}_{i=1}^{4g-4},x_0)$ as in Lemma \ref{S2lem}; with this choice, the analytic continuation of $\Psi$ has right monodromy $M_{\a_j}$ along $\a_j$ and $M_{\b_j}$ along $\b_j$. To show that these matrices  define an $SL(2,\C)$ representation of $\pi_1(\CC,x_0)$, one needs to verify the relation $\prod_{j=1}^g M^{-1}_{\b_j}M^{-1}_{\a_j} M_{\b_j}M_{\a_j} =\1$ implied by 
(\ref{relation}) (recall that monodromy matrices form an anti-representation of $\pi_1(\CC,x_0)$).

Now, the fundamental relation (\ref{relpunct}) in $\pi_1(\CC\setminus\{x_i\}_{i=1}^{4g-4},x_0)$ implies
\be
\prod_{j=1}^g M^{-1}_{\b_j}M^{-1}_{\a_j} M_{\b_j}M_{\a_j} \cdot \prod_{k=1}^{4g-4} M_{\g_k} = \1\;.
\ee
Since $M_{\g_k} = i\1$ and the number $ 4g-4$ of zeros  is a multiple of $4$,  the second product in this formula equals $\1$ and the  matrices $\{M_{\a_i}, M_{\b_i}\}$ define an $SL(2,\C)$ representation of 
$\pi_1(\CC,x_0)$. 
\QED

For later convenience, we  introduce the following combinations of the transition matrices and matrices $\sigma_j$:
\bea
T_{j}  ^{z_1,z_2}(z)= T(z_1,z) \sigma_j  T(z,z_2) = \Psi(z_1) \Lambda_j(z) \Psi^{-1}(z_2)\ ,\qquad j \in\{ 3,\pm\}.
\label{defTm}
\eea
The following notations will be also used:
\be
\l_{j}^{z_1,z_2}(z):=\tr T_j^{z_1,z_2}(z)\;,  \hskip0.7cm j\in\{3,\pm\}
\la{lj12def}
\ee
\be
\l_j^{z_1,z_2}:=\tr(\s_j T(z_1,z_2))= \l_{j}^{z_1,z_2}(z_1)=\l_{j}^{z_1,z_2}(z_2)\;,\hskip0.7cm j\in\{3,\pm\}
\la{lj12}
\ee
\be
\l_j^{\g}(z):=\l_j^{z_0+\P_\g,z_0}(z) =\tr[T(z_0+\P_\g ,z)\s_j T(z,z_0)]\;,\hskip0.7cm j\in\{3,\pm\}\;.
\la{ljgdef}
\ee

\section{The Poisson algebra of monodromies: Goldman bracket}
\label{SecPoi}
The goal of this section is to show that the homological symplectic form \eqref{fundam} on the space $\Q_g^0$ implies the 
 Goldman Poisson bracket for traces of monodromy matrices of the equation (\ref{matrix1}) on $\CV_g$. The monodromy map $\Q_g\to \CV_g$ is in fact a composition of two maps.
 The first one is the map  $\Q_g\to \Proj_g$ using the Bergman projective connection as the base, and the second one  is the map  $\Proj_g \to \CV_g$ 
 defined by the monodromy  of equation (\ref{matrix1}).

We start from reformulating the Poisson bracket (\ref{Poissonu}) as the following theorem:
\begin{theorem}
\label{thmPoissonu} 
The homological Poisson bracket of the potential $u(z)$ \eqref{uq} can be expressed as follows;
\be
\f{4\pi i}{3}\{u(z), u(\zeta)\}=\L_z h^{(\zeta)}(z)- \L_\zeta h^{(z)}(\zeta)\;,
\label{Poissonu1}
\ee
where the differential operator $\L_z$ is given by:
\be
\L_z:= \frac 1 2 \pa_z^3 - 2u(z) \pa_z  -u_z(z) \;.
\label{defLz}\ee
\end{theorem}
The differential operator $\L_z$ arises in the theory of KdV equation \cite{BaBeTa} where it is known as the Lenard operator and also appears in the theory of projective structures on Riemann surfaces (see   formula (9) in \cite{Gunning1}).

Using the definition of the 
Poisson bracket, for any closed loop $\g$ on $\CC$, the variational formula (\ref{varQ}) of $\bc(z)$ can be expressed in terms of $u(z)$ as follows: 
\be
\left\{u(z),\,\oint_{\g} v\right\}=\frac{3}{4\pi i} \oint_\g h(z,t)dt
\mathop{\equiv}^{\eqref{defHg}}
\frac{3}{4\pi i} \f{H^{(\g)}(z)}{v(z)}\;.
\label{uint1}
\ee 
Eq.  (\ref{uint1}) holds provided that the contour defining the coordinate $z(x)=\int_{x_1}^xv$ does not intersect $\gamma$.

The Poisson bracket between transition matrices of an arbitrary matrix
differential equation $\Psi_z=U(z)\Psi$ can be written as follows  (see \cite{FadTak} for explanation of the notation $\{\  \overset{\otimes}{,} \ \}$)
 \be
\left\{T(z_1,z_2)\overset{\otimes}{,}T(\zeta_1,\zeta_2)\right\}
=\int_{z_1}^{z_2}\int_{\zeta_1}^{\zeta_2}  T(z_1,z)\otimes T(\zeta_1,\zeta)\left\{U(z)\overset{\otimes}{,}U(\zeta)
\right\} T(z,z_2)\otimes T(\zeta,\zeta_2) \d\z \d z\;.
\label{Poisstrans}
\ee
Analogously, for an arbitrary   scalar function $f$ we have 
\be
\left\{T(z_1,z_2), f\right\}= \int_{z_1}^{z_2} T(z_1,z) \left\{U(z),f\right\}  T(z,z_2) \d z\;. 
\label{univ}\ee
In our case  
\be
U(z)= u(z)\sigma_-+\sigma_+
\label{Ustruct}
\qquad
\left\{U(z)\overset{\otimes}{,}U(\zeta)
\right\} = \{u(z), u(\zeta)\}\, \s_-\otimes\s_-
\ee
and substitution of    (\ref{uint1}) into (\ref{univ}) gives the following
\begin{lemma}
Assume the contour $[z_1,z_2]$ does not intersect $\gamma$ and $z_{1,2}$ remain constant. Then 
\be
\f{4\pi i}{3}\left\{T(z_1,z_2),\P_\g \right\}=\int_{z_1}^{z_2} \Tm^{z_1,z_2}(z)  H^{(\g)}(z)\d z 
\label{transhom}
\ee
where the differential $H^{(\g)}(z)$ is given by  (\ref{defHg})  and  $\Tm^{z_1,z_2}$ is defined by \eqref{defTm}.
\end{lemma}

\begin{corollary}
For any two non-intersecting closed contours $\g$ and $\gt$ on $C$, 
\be
\f{4\pi i}{3}\left\{\tr M_{\gt},\P_\g \right\}=\int_{\gt} H^{(\g)}(z) \tr [\Ld(z) M_{\gt}] \d z\;.
\label{trmon}
\ee
\end{corollary}

\begin{lemma}\label{auxl}
Let two non-intersecting arcs $l_1=[z_1,z_2]$ and $l_2=[\zeta_1,\zeta_2]$ lie entirely within the fundamental domain $C_0$; then 
\bea
\f{8 \pi i}{3} &\&\left\{T(z_1,z_2)\overset{\otimes}{,}T(\zeta_1,\zeta_2)\right\} = 
\cr
=&\&\int_{\zeta_1}^{\zeta_2} \left[ \Tm^{z_1,z_2}h_{zz}^{(\zeta)}+ \T3 ^{z_1,z_2} h_{z}^{(\zeta)}-2 
\left(\Tp^{z_1,z_2}+u(z)\Tm^{z_1,z_2}\right) h^{(\zeta)}\right]\Big|_{z=z_1}^{z=z_2}\otimes  \Tm^{\zeta_1,\zeta_2} \d\zeta - (z\leftrightarrow \zeta)\;.
\label{TTzz}
\eea

\end{lemma}

\noindent {\bf Proof.}
The expression (\ref{Poisstrans}) can be simplified using \eqref{Ustruct} 
\be
\left\{T(z_1,z_2)\overset{\otimes}{,}T(\zeta_1,\zeta_2)\right\}=
\int_{z_1}^{z_2}\int_{\zeta_1}^{\zeta_2}  \Tm^{z_1,z_2}(z) \otimes  \Tm^{\zeta_1,\zeta_2}(\zeta) \{u(z),u(\zeta)\}\, \d z \d\zeta
\label{Poisstrans1}
\ee
where $ \Tm^{z_1,z_2}(z)$ is defined in (\ref{defTm}).
As before the loops bounding the fundamental domain 
are chosen according to Lemma \ref{S2lem}.
Substituting (\ref{Poissonu1}) into (\ref{Poisstrans1}) we get
$$
\f{4\pi i}{3}\int_{l_1}\int_{l_2} \Tm^{z_1,z_2}(z) \otimes  \Tm^{\zeta_1,\zeta_2}(\zeta) \{u(z),u(\zeta)\}\,\d z \d\zeta
$$
\be
= 
\int_{l_2}\left(\int_{l_1}  \Tm^{z_1,z_2}\left( \L_z  h^{(\zeta)}\right)\d z\right)\otimes \Tm^{\zeta_1,\zeta_2} (\zeta) \d\zeta-
\int_{l_1} \Tm^{z_1,z_2}(z) \otimes\left(\int_{l_2} \Tm^{\zeta_1,\zeta_2}\left(\L_\zeta h^{(z)}\right) \d\zeta\right) \d z
\label{aux2}
\ee
where the order of integration may be interchanged because the paths $l_1$ and $l_2$ are non-intersecting.

The first integral in both of these double integrals can now be performed explicitly; namely,
 \be
\int_{l_1} \Tm^{z_1,z_2} \L_z h^{(\zeta)} \d z= \int_{z_1}^{z_2} \Tm^{z_1,z_2}(z) \left(\f{1}{2}h^{(\zeta)}_{zzz}-2 h^{(\zeta)}_{z}u(z) - h^{(\zeta)} u_z\right) \d z\;.
\label{aux1}\ee
Integrating the first term by parts three times we find
$$
\int_{z_1}^{z_2} h^{(\zeta)}_{zzz} \Tm^{z_1,z_2}\d z= -\int_{z_1}^{z_2} (\Tm^{z_1,z_2})_{zzz} h^{(\zeta)}\d z+ \left( h^{(\zeta)} (\Tm^{z_1,z_2})_{zz}
- h^{(\zeta)}_z  (\Tm^{z_1,z_2})_{z}+ h^{(\zeta)}_{zz}\Tm^{z_1,z_2} \right)\big|_{z_1}^{z_2},
$$
and integrating the second term by parts once,
 $$ \int_{z_1}^{z_2} h^{(\zeta)}_{z} u\, \Tm^{z_1,z_2} \d z= -\int_{z_1}^{z_2} h^{(\zeta)}(u \,\Tm^{z_1,z_2}(z) )_z \d z+ h^{(\zeta)}_{z}u \,\Tm^{z_1,z_2}(z)\big|_{z_1}^{z_2} \;.$$  
After these integrations by parts the integral term in  (\ref{aux1})   vanishes by the third order 
equation (\ref{3rdor}) for $\Ld$. The remaining boundary terms in (\ref{aux1}) equal  
$$
\f{1}{2}\left( \Tm^{z_1,z_2}  h^{(\zeta)}_{zz}+ \T3 ^{z_1,z_2}  h^{(\zeta)}_z - 2(\Tp^{z_1,z_2} + u \Tm^{z_1,z_2}) h^{(\zeta)}(z)\right)\big|_{z_1}^{z_2}
$$
which is exactly the first term in the r.h.s. of (\ref{TTzz}). Similarly the second term in the r.h.s of  
(\ref{TTzz}) is the result of performing the integration over $\zeta$ 
in the second double integral of the formula (\ref{aux2}).
\QED

\subsection{Brackets between traces of monodromy matrices}
In this section we show that the Poisson bracket (\ref{Poissonu1}) for the potential $u$ implies the Goldman Poisson bracket for monodromy matrices. The main result is Thm. \ref{sympGold}, to which we are going to arrive  via Thm. \ref{nonint} and Thm. \ref{intth}. 

\begin{lemma}
Let $\gamma$ be any closed loop on $\CC$, and $[\zeta_1,\zeta_2]$ be an arc not intersecting $\gamma$, then
\bea
\f{8 \pi i}{3}&\& \left\{\tr \,M_\gamma\,,\, \tr\, T(\zeta_1,\zeta_2)\right\}
\cr
=&\&-\l_-^{\zeta_1,\zeta_2}\int_{z\in\g} \l_-^\g(z) \left(\h_\zeta(\zeta,z)\big|_{\zeta=\zeta_1}^{\zeta=\zeta_2}\right)\d z-
\l_3^{\zeta_1,\zeta_2}\int_{z\in\g} \l_-^\g(z) \left(\h(\zeta,z)\big|_{\zeta=\zeta_1}^{\zeta=\zeta_2}\right)\d z
\cr 
&\& +2\l_+^{\zeta_1,\zeta_2}\int_{z\in\g} \l_-^\g(z)\int_{\zeta_1}^{\zeta_2} \h(\zeta,z) \d\zeta \d z
\cr
&\& +2\l_-^{\zeta_1,\zeta_2}\left[u(\zeta_2)\int_\g \l_-^\g(z)\left(\int_{x_1}^{\zeta_2} \h(z,\zeta)\d\zeta\right)\d z
-u(\zeta_1)\int_\g \l_-^\g(z)\left(\int_{x_1}^{\zeta_1} \h(z,\zeta)\d\zeta\right)\d z\right]\;.
\label{MgT}
\eea
\end{lemma}
\noindent {\bf Proof.}  
Taking the trace  of  both factors in  the tensor product appearing in  (\ref{TTzz}) we obtain the expression
\bea
\f{8 \pi i}{3}&\& \left\{\tr \,T(z_1,z_2){,}\, \tr\, T(\zeta_1,\zeta_2)\right\} =
\cr
=&\&\l_-^{z_1,z_2}(z_1)\int_{\zeta_1}^{\zeta_2}\l_-^{\zeta_1,\zeta_2}(\zeta)
\left(h_z(z,\zeta)\big|_{z=z_1}^{z=z_2}\right) \d\zeta+
\l_3^{z_1,z_2}(z_1)\int_{\zeta_1}^{\zeta_2}\l_-^{\zeta_1,\zeta_2}(\zeta)
\left(h(z,\zeta)\big|_{z=z_1}^{z=z_2}\right) \d\zeta
\cr
&\&-2\l_+^{z_1,z_2}(z_1)\int_{\zeta_1}^{\zeta_2}\l_-^{\zeta_1,\zeta_2}(\zeta)\left(\int_{z_1}^{z_2}h(z,\zeta) \d z \right)\d\zeta
\cr
&\&
-2\l_-^{z_1,z_2}(z_1)\int_{\zeta_1}^{\zeta_2}\l_-^{\zeta_1,\zeta_2}(\zeta)\left[\left(u(z)\int_{x_1}^z h(z',\zeta)\d z'\right)\big|_{z=z_1}^{z=z_2}\right]\d\zeta-
(z\leftrightarrow\zeta)\;.
\label{llzz}
\eea
 The symbol $ (z\leftrightarrow\zeta)$ means that all $z,z', z_1,z_2$ are interchanged with $\z,\z', \z_1,\z_2$, respectively.
 This expression   should be evaluated at $z_2=z_0, \ \ z_1  = z(x_0^\gamma) = z_0 + \P_\g$. 
 Thus there appears an additional implicit dependence on the homological coordinate $\P_\g$ which needs to be taken into account.
Namely, 
$
M_\g= T(z_0+\P_\g,z_0)\;,
$
and 
\be
\label{llzz1}
\{\tr M_\g\,,\,\tr T(\zeta_1,\zeta_2)\} =\le\{\tr T(z_1,z_2)\,,\,\tr T(\zeta_1,\zeta_2)\ri\}\big|_{
z_1=z_0+\P_\g
\atop 
z_2 = z_0}
+\tr T_{z_1}(z_0+\P_\g,z_0)\{\P_\g,\tr T(\zeta_1,\zeta_2)\}
\ee
where $T_{z_1}$ denotes the derivative with respect to the first argument.
The expression for $\{\P_\g,\tr T(\zeta_1,\zeta_2)\}$ is given by (\ref{transhom}) while
$$
\tr T_{z_1}(z_0+\P_\g,z_0)= \tr\le(\Psi' (z_1) \Psi^{-1}(z_0)\ri) = \tr\le(\Psi (z_1) (\s_+ + u(z_1)\s_-)\Psi^{-1}(z_0)\ri)   
$$
\be
= u(z_0)\l_-^\g(z_0)+\l_+^\g(z_0)
\label{trTprime}
\ee
using the fact that $u(z)$ is single-valued on $\CC$.

The first two integrals in the r.h.s. of (\ref{llzz}) vanish when the arc $[z_1,z_2]$ closes because $\h(z,\zeta)$ is a single--valued function on $\CC\times \CC$.
The third and fourth integrals together give
$$
-2 (\l_+^\g(z_0)+ u(z_0) \g_-^\g(z_0))\int_{\zeta_1}^{\zeta_2}\l_-^{\zeta_1\zeta_2}(\zeta)\left(\int_{z_0+\P_\g}^{z_0} h(\zeta,z) \d z\right) \d\zeta
$$
\be
=2T_{z_1}(z_0+\P_\g,z_0)\int_{\zeta_1}^{\zeta_2}\l_-^{\zeta_1\zeta_2}(\zeta)H^{(\gamma)}(\zeta) \d\zeta
\label{addterm}\ee
which, by \eqref{transhom}, cancels with the second term in the right side of \eqref{llzz1}.
 Thus the r.h.s of (\ref{MgT}) consists only of the term indicated by  $(\zeta\leftrightarrow z)$ in (\ref{llzz}).
 \QED


\begin{theorem}\label{nonint}
If $\g$ and $\gt$ are non-intersecting loops on $\CC$, the traces of the corresponding monodromy matrices Poisson-commute:
\be
\{{\rm tr} M_{\g},{\rm tr}  M_{\gt}\}=0\;.
\ee
\end{theorem}
{\bf Proof.} The trace of a monodromy matrix does not depend on the basepoint, and the arc $[\zeta_1,\zeta_2]$ in (\ref{MgT}) can be closed to form a loop by assuming $\zeta_2=\zeta_0$,
 $\zeta_1=\zeta_0+\P_{\gt}$.
Therefore
\be
{\f{8\pi i}{3}}
 \{{\rm tr} M_{\g},{\rm tr}  M_{\gt}\}=
{ \f{8\pi i}{3}}
  \{{\rm tr} M_{\g},\tr T(\zeta_1,\zeta_2)\}\big|_{\zeta_1=\zeta_0+\P_\gt
  \atop \zeta_2=\zeta_0}
+
{\f{8\pi i}{3}}
\{{\rm tr} M_{\g},\, \P_\gt\} \tr T_{z_1}(\zeta_0+\P_\gt,\zeta_0)
\label{ggt}
\ee
where the first term is given by (\ref{MgT}) with the substitution $\zeta_1=\zeta_0+\P_\gt$, $\zeta_2=\zeta_0$.
The first two integrals in the r.h.s. of (\ref{MgT}) vanish when $\zeta_1=\zeta_0+\P_\gt$ and $\zeta_2=\zeta_0$ because
$h_\zeta(z,\zeta)\big|_{\zeta_0+\P_\gt}^{\zeta_0}=h(z,\zeta)\big|_{\zeta_0+\P_\gt}^{\zeta_0}=0$.

Denote $ \displaystyle \l_i^\gt=\l_i^{\zeta_0+\P_\gt,\zeta_0}\;,$
 $i=(3,\pm)$, and recall that $[\zeta_0+\P_\gt,\zeta_0]=-\gt$. The remaining terms in (\ref{MgT}) are
\be
2(\l_+^\gt+u(\zeta_0)\l_-^\gt)\int_{z\in\g}\l_-^\g (z) H^{(-\gt)}(z)\d z
\label{ggtt}
\ee
 which (taking into account that $H^{(-\gt)}=-H^{(\gt)}$) cancel against the second term in (\ref{ggt})  due to \eqref{trTprime} and (\ref{transhom}).  \QED

\begin{theorem}\label{intth}
Let $\g$ and $\gt$ be two closed contours on $\CC$ intersecting transversally at one point  with $\gt\circ\g=1$. Denoting the monodromy matrices
of equation (\ref{matrix1}) corresponding to the same basepoint $z_0$ by $M_\g$ and $M_\gt$, we have
\be
\{\tr M_\g,\,\tr M_\gt\}=\f{1}{2}(\tr M_\g M_\gt-\tr M_\g M^{-1}_\gt)\;.
\label{Gold}
\ee
\end{theorem}

\begin{figure}
\centering
\resizebox{0.3\textwidth}{!}{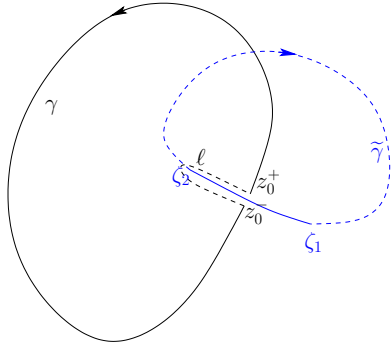}
\captionof{figure}{The closed loop $\g$ intersecting the interval $[\zeta_1,\zeta_2]$ and the arc $\ell$.
 }
\label{loopinterval}
\end{figure}

Before proceeding to the proof  of the theorem, we  prove  the following auxiliary statement:
\begin{lemma}\label{intl}
Let $[\z_1,\z_2]$ be an arc in the fundamental domain and $z_0$ a point on it. Let $\ell$ be the arc connecting points $z_0^-$ and $z_0^+$ on the opposite sides of $[\z_1,\z_2]$ at $z_0$ as indicated in Fig. \ref{loopinterval}. Then
the following formulas hold:
\be
\int_\ell \l_-^\g(z) h(z,\zeta_2)\d z =-\f{4\pi i}{3} \l_3^\g(\zeta_2)\;,
\label{intlh}
\ee
\be
\int_\ell \l_-^\g(z) h_\zeta(z,\zeta_2)\d z =-\f{8\pi i}{3} (u(\zeta_2) \l_-^\g(\zeta_2)- \l_+^\g(\zeta_2))\;,
\label{intlhz}
\ee
\be
\int_\ell \l_-^\g(z)\left(\int_{\zeta_1}^{\zeta_2} h(z,\zeta)d\zeta\right)\d z =\f{4\pi i}{3}\l_-^\g(\zeta_2)+\f{2\pi i}{3} (u(z_0) \l_-^\g(z_0)+ \l_+^\g(z_0))\;.
\label{intintlh}
\ee

\end{lemma}
\noindent {\bf Proof.} Recall  the asymptotics \eqref{asymph} of $h(z,\zeta)$ as $\z\to\zeta$:
\be
h(z,\zeta)=\f{1}{(z-\zeta)^4}-\f{2}{3}\f{u(\zeta)+1}{(z-\zeta)^2}+\f{u_\zeta(\zeta)}{3(z-\zeta)} +\dots
\label{hdia}
\ee

The integrand in (\ref{intlh}) and in  (\ref{intlhz}) do not have a jump at $z_0$. Thus the integral (\ref{intlh}) along the closed contour $\ell$ is given by the residue at $z=\zeta_2$:
\be
\res\big|_{\zeta_2} \l_-^\g(z) h_\zeta(z,\zeta_2)dz=\f{1}{6}(\l_-^{\g}) ''' (\zeta_2)-\f{2}{3}(u(\zeta_2)+1) (\l_-^{\g})'(\zeta_2) +\f{1}{3}u_\zeta(\zeta_2) \l_-^\g(\zeta_2)
\ee
which by the third order equation for traces (\ref{3rdor}) is $-\f{2}{3} (\l_-^\g)'(\zeta_2)=\f{2}{3} (\l_3^\g)'(\zeta_2)$.
The contour $\ell $ is negatively oriented, and therefore one multiplies by $-2\pi i$ to get (\ref{intlh}).

The integral (\ref{intlhz}) is computed by differentiation of (\ref{intlh}) with respect to $\zeta_2$ and by   using (\ref{derF}), (\ref{derFtt}).

Finally, to compute the integral (\ref{intintlh}), we first interchange the order of integration over intersecting contours $\ell$ and $[\zeta_1,\zeta_2]$. 
Changing   the order of integration in a singular integral gives rise to an additional term which can be computed using the universal formula which is valid
for arbitrary  arcs $s_1$ and $s_2$ intersecting at the point $x_0$ with intersection index $s_1\circ s_2=1$:
\be
\left(\int_{s_1}\d x\int_{s_2}\d y -\int_{s_2}\d y\int_{s_1}\d x\right)\f{f(x)g(y)}{(x-y)^n} =\f{2\pi i}{(n-1)!} \left(\sum_{k=0}^{n-2}\p_x^k\p_y^{n-2-k}\right) f(x)g(y)\big|_{x=y=x_0}
\label{interchange}
\ee
where $f(x), g(x)$ are locally analytic expressions at the intersection point.

The additional terms arising from interchanging the order of integration in (\ref{intintlh}) are computed by applying (\ref{interchange}) to the asymptotic expansion (\ref{hdia}) and re-expressing using (\ref{derFttt}):
$$
\int_\ell \l_-^\g(z)\left(\int_{\zeta_1}^{\zeta_2} h(z,\zeta)\d\zeta\right)\d z=\int_{\zeta_1}^{\zeta_2}\left(\int_\ell  \l_-^\g(z) h(z,\zeta)\d z\right)\d\zeta
$$
\be
+\f{2\pi i}{3}(u(z_0)\l_-^\g(z_0) +\l_+^g(z_0)+2\l_-^\g(z_0))\;.
\label{iiint}
\ee

The integral in the r.h.s. of (\ref{iiint}) can be computed explicitly to give
$$
\int_l \l_-^\g(z) h(z,\zeta)\d z=-\f{4\pi i}{3} \l_3^\g(\zeta)
$$
for $\zeta$ lying between $z_0$ and $\zeta_2$, and $0$ if $\zeta$ lies between $\zeta_1$ and $z_0$.

Therefore,
\be
\int_{\zeta_1}^{\zeta_2}\left(\int_l \l_-^\g(z) h(z,\zeta)\d z\right)\d\zeta=-\f{4\pi i}{3} (\l_-^\g(\zeta_2)-\l_-^\g(z_0))
\label{int2}
\ee
Combining (\ref{iiint}) with  (\ref{int2}), one gets (\ref{intintlh}).
\QED

\noindent {\bf Proof of Thm. \ref{intth}.}  The right-hand side of (\ref{Gold}) can be rewritten as follows:
\be
\f{1}{2}\Big(\tr M_\g M_\gt-\tr M_\g M^{-1}_\gt\Big )= \l_+^\g\l_-^\gt +\l_-^\g\l_+^\gt +\f{1}{2}\l_3^\g\l_3^\gt
\label{Gold1}
\ee
where, as before, $\l_i^\g=\tr (M_\g\s_i)$ and  $\l_i^\gt=\tr (M_\gt\s_i)$.
Although the individual terms $\l_i^\g$ in (\ref{Gold1}) are not well-defined functions on the character variety $\CV_g$, as they depend on the point $z_0$, the combination in the r.h.s. of \eqref{Gold1} is independent of $z_0$ because it can be expressed in terms of traces of monodromies.

Denote points of $\g$ lying on different sides of the arc $[\zeta_1,\zeta_2]$ by $z_0^\pm$ i.e. $\g=[z_0^+,z_0^-]$ (Fig. \ref{loopinterval}),
and introduce the closed contour $\ell$ which goes around $\zeta_2$ from $z_0^-$ to $z_0^+$  as shown in Fig. \ref{loopinterval}. 

Let  $[\zeta_1,\zeta_2]$ be  an arc which intersects the contour $\g$ at a point $z_0$.
In order to derive an analog of the formula (\ref{MgT}) for 
$\{\tr M_\g,\tr T(\zeta_1,\zeta_2)\}$, we substitute the contour $\g$ with the contour $\g\cup \ell $ which does not intersect
 $[\zeta_1,\zeta_2]$, noting $\tr M_\g=\tr M_{\g \ell}$. Then we can compute this bracket using the formula for non-intersecting paths
(\ref{MgT}). 

We compute the contribution of $\ell$ with Lemma \ref{intl}.
In a similar manner two other integrals can be computed. 
Choose a path from $x_1$ to $\zeta_1$ which does not intersect $\g$, then
$$
\int_\ell \l_-^\g(z)\left(\int_{x_1}^{\zeta_1} h(z,\zeta)\d\zeta\right)\d z=0
$$
where interchanging the order of integration does not result in any extra terms because the singularity of $h(z,\zeta)$ lies outside of $l$.
Analogously, the integral 
$$
\int_\ell\l_-^\g(z)\left(\int_{x_1}^{\zeta_2} h(z,\zeta)\d\zeta\right)\d z
$$
coincides with expression  (\ref{intintlh}).
Therefore, 
\bea
\f{8\pi i}{3}\{\tr M_\g, \tr T(\zeta_1,\zeta_2)\}
=&\& \{ \hbox{ expression (\ref{MgT}) where  $z$ is integrated along $\g$ between $z_0^+$ and $z_0^-$}  \}\cr
&\& + \{\hbox{the same expression where integration over $z$ goes along $l$}\}
\nonumber
\eea

The contribution of the contour $\ell $ is computed using Lemma \ref{intl}:
$$
\f{8\pi i}{3}\bigg(\l_-^{\zeta_1,\zeta_2} \l_+^\g (\zeta_2)+\l_+^{\zeta_1,\zeta_2} \l_-^\g (\zeta_2)+\f{1}{2}\l_3^{\zeta_1,\zeta_2} \l_3^\g (\zeta_2)\bigg)
+\f{4\pi i}{3}\bigg( u(\zeta_2)   \l_-^{\zeta_1,\zeta_2}  +\l_+^{\zeta_1,\zeta_2}\bigg )\bigg( u(z_0)   \l_-^\g(z_0)  +\l_+^\g(z_0)\bigg )\;.
$$
We now close the contour $[\zeta_1,\zeta_2]$ by setting $\zeta_1=\zeta_2+\P_\gt$ to get
$$
\{\tr M_\g, \tr M_\gt\}= \Big\{\tr T(z_0+\P_\g,z_0),\tr T(\zeta_2+\P_\gt,\zeta_1)\Big\}
$$
$$
=\Big\{\tr T(z_0+\P_\g,z_0) , \tr T(\zeta_1,\zeta_2)\Big\}\big|_{\zeta_1=\zeta_2+\P_\gt}+ 
\Big\{\tr T(z_0+\P_\g,z_0), \P_\gt \Big\}\;\tr T_1'(\zeta_2+\P_\gt,\zeta_1)\;.
$$
Since  $\g$ and $\gt$ intersect at one point with negative orientation in $\CC$, the intersection index of their images in $H_-$ under the homomorphism ${\bf g}$ 
(\ref{homom1})  equals $-1/2$  (due to normalization (\ref{intH-}) and the definition (\ref{abm1}) of the lift  of $H_1(\CC)$ to $H_-$). Therefore,  (\ref{PSfund}) implies 
$$
 \{\P_\g,\P_\gt\}=- \frac 1 2 \;,
$$
$$
\Big\{\tr T(z_0+\P_\g,z_0), \P_\gt\Big\}=\f{1}{2}\tr T_1'(z_0+\P_\g,z_0) +\f{3}{4\pi i} \int_{z\in \g} \l_-^\g (z) H^{(\gt)}(z)\d z\;.
$$
The following holds in analogy with the proof of Theorem \ref{nonint}.
$$
\Big\{\tr T(z_0+\P_\g,z_0) , \tr T(\zeta_1,\zeta_2)\Big\}\big|_{\zeta_1=\zeta_2+\P_\gt}+\f{3}{4\pi i} \tr T_1'(\zeta_2+\P_\gt,\zeta_2) \int_{z\in \g} \l_-^\g (z) H^{(\gt)}(z)\d z =0\;.
$$
\noindent
The remaining terms are 
$$
\{\tr M_\g,\; \tr M_\gt\} =
$$
$$
 =\f{1}{2}\left(u(\zeta_2)\l_-^{(\zeta_1+\P_\gt,\zeta_2)}+ \l_+^{(\zeta_1+\P_\gt,\zeta_2)}\right)
\Big(u(z_0)\l_-^\g(z_0)+\l_+^\g(z_0)\Big)
$$
$$
-\f{1}{2}T_1'(z_0+\P_\g,z_0)T_1'(\zeta_2+\P_\gt,\zeta_2)
$$
$$
+\left(\l_-^\gt(\zeta_2)\l_+^\g(\zeta_2)+\l_+^\gt(\zeta_2)\l_-^\g(\zeta_2)+\f{1}{2}\l_3^\gt(\zeta_2)\l_3^\g(\zeta_2)\right)
$$
where the first two terms cancel each other since
$$
T_1'(z_0+\P_\g,z_0)= u(z_0)\l_-^\g(z_0)+\l_+^\g(z_0) +T_1'(\zeta_2+\P_\gt,\zeta_2)
$$
$$
=u(\zeta_2)\l_-^{(\zeta_1+\P_\gt,\zeta_2)}+ \l_+^{(\zeta_1+\P_\gt,\zeta_2)}\;.
$$
The final result is
$$
\{\tr M_\g, \tr M_\gt\}=\l_-^\gt(\zeta_2)\l_+^\g(\zeta_2)+\l_+^\gt(\zeta_2)\l_-^\g(\zeta_2)+\f{1}{2}\l_3^\gt(\zeta_2)\l_3^\g(\zeta_2)
$$
which is the Goldman bracket corresponding to the initial point $\zeta_2$; the expression in the r.h.s. is independent of $\zeta_2$ and coincides with the Goldman bracket for paths intersecting at one point.
\QED
\paragraph{Bracket between traces of monodromies along two arbitrary loops.}
The main technical result of the paper now follows from the previous computation.
\begin{theorem}
\label{sympGold}
Let $\g,\gt\in \pi_1(\CC,x_0)$ and let $M: \pi_1(\CC,x_0) \to SL(2,\C)$ be  the monodromy (anti-)representation of the equation \eqref{matrix1} defined in Section \ref{secmon}.
The homological symplectic structure  \eqref{fundam}, or, equivalently, the canonical Poisson structure on $T^*\Mc_g$, implies the following Poisson bracket between traces of monodromy matrices:
\be
\{\tr M_\g,\;\tr M_\gt\}= \f{1}{2}\sum_{p\in \g\cap \gt} (\tr M_{\g_p\gt}-\tr M_{\g_p\gt^{-1}})
\label{Goldman}
\ee  
where $\g_p\gt$ and $\g_p\gt^{-1}$ are two ways to resolve the intersection point $p$ to get two new contours  $\g_p\gt$ and $\g_p\gt^{-1}$ for each $p$.
\end{theorem}

{\bf Proof.}  One can construct a sufficiently large set of functions which generically define a point of $\CV_g$ as follows. Choose a standard set of generators 
$\{\a_i,\b_i\}$ of $\pi_1(\CC,x_0)$  satisfying (\ref{relation}) and consider the following set of $g^2+2g$ loops:
\be
\SS=\{\a_i,\ \b_i,\; \ 1 \leq i \leq g\;; \;\;\; \a_i\a_j, \; 1\leq i<j\leq g\;;\; \;\; \a_i\b_j, \; 1 \leq i \leq j \leq g\}\;.
\label{test}
\ee

Any two loops from this set either do not intersect or intersect at one point. Therefore, (\ref{Goldman}) holds for any pair of loops 
$\g,\gt\in \SS $.
 The number  $g^2+2g$ of these loops is always greater than the
 number or functions ($6g-6$) required to verify (\ref{Goldman}) for any pair of loops. 

It remains to verify that generically differentials of $\tr M_\g$ for $\g\in \SS $ generate  $T^*\CV_g$, or in other words, knowing $\tr M_\g$ for
$\g\in \SS $, one should generically be able to reproduce the point of the character variety $\CV_g$ up to a finite choice.

A classical theorem of Vogt and Fricke (see \cite{Gold09}) for details) states that a pair of $SL(2,\C)$ 
generic matrices $(X,Y)$ is completely defined, up to a simultaneous conjugation, by three traces $(\tr X,\tr Y,\tr XY)$ as long as 
$\tr (X)^2 + \tr(Y)^2 + \tr(XY)^2 - \tr(X)\tr (Y)\tr(XY)\neq 4$.
The space $\CV_g$ is parametrized by the ``complex Fricke'' coordinates, which are $6g-6$ independent matrix entries of $M_{\a_i},\, M_{\b_i}$ for
$i=1,\dots,g-1$. To factor out the freedom of simultaneous conjugation of all monodromy matrices, one assumes that $M_{\a_g}$ is diagonal, and that $1$ is a 
fixed-point of $M_{\b_g}$. Under these assumptions fixing
 $M_{\a_i},\, M_{\b_i}$ for $i=1,\dots,g-1$ determines the matrices $M_{\a_g}$ and $M_{\b_g}$ by the relation (\ref{relation}).

We have used the fact that to fix an $SL(2,\C)$ matrix up to a sign (discrete freedom is not important for computing the local rank of the map), 
it is sufficient to know the trace of the products of this matrix with
three other generic matrices. Representing the matrix $A$ as $aI+b\s_++c\s_-+d\s_3$ and fixing
$\tr (A A_i)$ for $i=1,2,3$, one gets a system of three linear equations on $a,b,c,d$ and one quadratic equation ${\rm det} A=1$.

Let us show now that all Fricke coordinates can be determined in general by knowing $\tr  M_\g$, $\g\in \SS $.
We start with the matrices $M_{\a_g}$ and $M_{\b_g}$ which are assumed to be diagonal and have fixed point equal to $1$ respectively. The set (\ref{test}) contains a triple $(\tr M_{\a_g},\; \tr M_{\b_g},\; \tr M_{\a_g} M_{\b_g})$, it is enough to reproduce $M_{\a_g}$ and $M_{\b_g}$.
Now to reproduce the Fricke coordinates contained in $M_{\a_i}$ for $i=1,\dots,g-1$, use the following subset of $\SS $: $(\a_i, \a_i\a_g, \a_i\b_g)$;
the matrices $M_{\a_g}$ and $M_{\b_g}$ are already fixed and knowing the triple $(\tr M_{\a_i},\; \tr M_{\a_i}M_{ a_g},\; \tr M_{\a_i} M_{\b_g})$, one can generically determine 
$M_{\a_i}$ up to a sign.

Therefore all the matrices $M_{\a_i}$, $i=1,\dots,g$ and $M_{\b_g}$ have been determined, and it remains to determine $M_{\b_j}$ for $j=1,\dots,g-1$. For  each $j$ the set (\ref{test}) 
contains $\b_j$ and $\a_i\b_j$ for $i=1,\dots,j$, and therefore, for each $j=2,\dots,g-1$ one knows $\tr M_{\b_j}$ and at least two other traces of products of $M_{\b_j}$
with known matrices; thus the matrices $ M_{\b_j}$ are determined up to a binary choice.

The only remaining matrix is $M_{\b_1}$ for which one knows $\tr M_{\b_1}$ and $\tr M_{a_1} M_{\b_1}$. In addition, the relation (\ref{relation})
involves matrices which have all been determined with the exception of $M_{\b_1}$, and thus it can be represented in the form $A M_{\b_1}B=   M_{\b_1}$ where $A$ and $B$ are 
two known matrices in $SL(2, \C)$. Generically, this equation determines $M_{\b_1}$ up to rescaling where the multiplicative constant can be found from $\tr M_{\b_1}$ or $\tr M_{\a_1} M_{\b_1}$.

Therefore, in an open neighbourhood of $\CV_g$ the values of $\{\tr M_\g\}$ for $\g\in \SS $ determine the point of the character variety $\CV_g$ up to a finite 
choice. On this neighbourhood our Poisson bracket coincides with the Goldman's bracket  $\{\cdot,\cdot\}_G$, and this relation can be extended to the whole $\CV_g$ by analyticity.
\QED

\subsection{Comparison with the results of Kawai }

\label{SectBers}
In Kawai's paper \cite{Kawai} the reference projective connection is chosen as the Bers  connection $\proj _{Bers}^{\CC_0,\eta}$, discussed in Section \ref{Bersproj} ($\CC_0$ is a chosen point in the Teichm\"uller space). 
Note that  the Bers  connection is  non-holomorphic with respect to the
moduli of the "initial" Riemann surface $\CC_0$, but it is holomorphic with respect to the moduli of the quadratic differential $\eta$ .
The main result of \cite{Kawai} was the statement that the symplectic structure induced on $\Proj_g$ from the canonical symplectic structure on $T^*\Mc_g$ by choosing  Bers'
as the reference projective connection also implies the Goldman bracket on the character variety $\CV_g$ under the monodromy map from $\Proj_g$ to $\CV_g$.

Comparing this theorem to our Theorem \ref{sympGold}, the symplectic structures induced on $\Proj_g$  from the canonical symplectic structure on 
$T^* \Mc_g$ by the choice of the   reference Bergman and Bers projective connections are equivalent. Therefore, there exists a generating function $G_{Berg}^{Bers,\,\CC_0}$
(which for any point $\CC_0\in \Te_g$ is a holomorphic function on the Teichm\"uller space $\Te_g$) of the corresponding change of Lagrangian embedding.
Thus, in analogy to Sections \ref{BCdm} and \ref{Lagrchange}, the following statement holds:
\begin{corollary}
There exists a function $G_{Berg}^{Bers}$ on the Teichm\"uller space such that 
\be
\delta_\mu G_{Berg}^{Bers}=\f{1}{2} \int_{\CC}(S_B-S_{Bers}^{\CC_0,\eta}) \mu 
\la{BERBER}
\ee
where $\mu$ is an arbitrary Beltrami differential on $\CC$. 
\end{corollary}

Equations of the type (\ref{BERBER}) arise in the theory of holomorphic factorization of the determinant of the Laplacian using a quasi-Fuchsian analog of Selberg's 
zeta-function \cite{MakTeo} (a natural analog of the Bowen-Zograf $F$-function (\ref{defF}) in the quasi-Fuchsian case).
To the best of our knowledge the solution of   (\ref{BERBER}) is not known. However, on the basis of the results of \cite{MakTeo},
\cite{Kim2005} and \cite{DhoPho1986}, we  propose the following conjectural expression for the $G_{Berg}^{Bers}$. 
The Selberg zeta-function
corresponding to the quasi-Fuchsian group $\Gamma_{\CC_0,\eta}$ can be defined in analogy to the ordinary Selberg zeta-function  by  (see \cite{MakTeo} for details), 
\be
{\mathcal Z}[\Gamma_{\CC_0,\eta}](s) = \prod_{{\gamma}}\prod_{m=0}^\infty (1-q_\gamma^{s+m})
\la{Selberg}
\ee
where $\g$ runs over all distinct primitive classes in $\Gamma_{\CC_0,\eta}$ excluding the identity; $q_\gamma$ is the multiplier of the group element $\g$.
The expression (\ref{Selberg}) is defined where it converges, and it is assumed to be analytically extended to the maximal domain within the 
space of quasi-Fuchsian groups.

The fundamental domain of the group  $\Gamma_{\CC_0,\eta}$ in $\C$ 
can be represented as a disjoint union of two simply-connected fundamental domains. One is the fundamental domain of the Riemann surface $\CC=\CC_0^\eta$,
and the other is the fundamental domain of the Riemann surface $\overline{\CC_0}$ (which is the mirror image of $\CC_0$, see \cite{Bers59}). Denote their period matrices (corresponding to the Torelli marking defining the Bergman projective connection on $\CC$) by $\O$ and $\tilde{\O}_0=-\overline{\O}_0$.

\begin{conjecture}
The solution of (\ref{BERBER}) is given by
\be
G_{Berg}^{Bers} =- 6\pi i\log\frac{{\mathcal Z}'[\Gamma_{\CC_0,\eta}](1)}{\det (\O-\overline{\O}_0)}
\la{conjfor}
\ee
\end{conjecture}

The main motivation for proposing expression (\ref{conjfor}) comes from the consideration of the case $\eta=0$. 
Namely, according to the formula of 
d'Hoker-Phong \cite{DhoPho1986} that relates the Selberg zeta-function to the determinant of the Laplace operator in hyperbolic metric on $\CC$, and the variational formula
for the determinant of Laplacian \cite{TakZog}, the (real-valued) solution of the equation 
$
\delta_\mu G=\f{1}{2} \int_{\CC_0}(S_B- S_F)\mu\;
$
is obtained in terms of the the usual Selberg zeta-function $\mathcal Z$ corresponding to the Fuchsian group of $\CC_0$ i.e.
\be
\delta_\mu \left(6\pi i \,\log\f{{\mathcal Z}'(1)}{{\rm det}\Im \Omega_0}\right)=\f{1}{2} \int_{\CC_0}(S_F- S_B)\mu\;.
\la{dmuZ}
\ee

The left-hand side of (\ref{dmuZ}) is the limit of the right-hand side of (\ref{conjfor}) as $\eta=0$ while the integral in the right-hand side of (\ref{conjfor}), in the limit as $\eta\to 0$, is the integral in the right-hand side of (\ref{dmuZ}); therefore the conjectured expression (\ref{conjfor}) indeed satisfies the required equations for $\eta=0$, i.e. in the "Fuchsian" case.
In the case of quasi-Fuchsian groups, the solution of 
(\ref{BERBER}) can not be real because it must be holomorphic with respect to the moduli of $\eta$ (but not with respect to the moduli of $\CC_0$!). The expression (\ref{conjfor}) is, in our opinion, a natural way 
to analytically continue the  expression $6\pi i\, \log {\mathcal Z}'(1)/({\rm det}\Im \Omega)$ to the quasi-Fuchsian case.

\end{document}